\documentclass[12pt,a4paper]{amsart}
\usepackage{geometry}
\usepackage{tabls}
\usepackage{url}
\usepackage[utf8]{inputenc}
\usepackage{amsfonts}
\usepackage{amssymb}
\usepackage{mathrsfs}
\usepackage{amsmath}
\usepackage{amsthm}
\usepackage{amscd}
\usepackage{fancyhdr}
\usepackage{enumerate}
\usepackage{paralist}
\usepackage{graphicx}
\usepackage{cite}

\usepackage{footmisc}
\numberwithin{equation}{section}

\theoremstyle{plain}

\newtheorem{Thm}[equation]{Theorem}
\newtheorem{lem}[equation]{Lemma}

\newtheorem{rem}[equation]{Remark}

\newtheorem{conj}[equation]{Conjecture}

\begin{document}

\title{Diophantine approximation  of multiple zeta-star values}

\author{Jiangtao Li}

\email{lijiangtao@csu.edu.cn}
\address{Jiangtao Li \\ School of Mathematics and Statistics, HNP-LAMA, Central South University, Hunan Province, China}

\begin{abstract}
     The set of multiple zeta-star values is a countable dense subset of the half line $(1,+\infty)$. In this paper, we establish some classical Diophantine type results for the set of multiple zeta-star values. Firstly, we give a criterion to determine whether a number is a multiple zeta-star value. Secondly, we establish the zero-one law for the set of multiple zeta-star value. Lastly, we propose a conjecture for the set of multiple zeta-star values, which strengthens the original  zero-one law.
     
\end{abstract}

\let\thefootnote\relax\footnotetext{
2020 $\mathnormal{Mathematics} \;\mathnormal{Subject}\;\mathnormal{Classification}$. 11M32, 11J83.\\
$\mathnormal{Keywords:}$  Multiple zeta-star values, Diophantine approximation. }

\maketitle

\section{Introduction}\label{int}
Multiple zeta values are defined by 
 \[
      \zeta(k_1,\cdots,k_r)=\sum_{n_1> \cdots> n_r\geq 1}\frac{1}{n_1^{k_1}\cdots n_r^{k_r}}, k_1\geq 2,k_2,\cdots, k_{r}\geq 1.
      \]
      For a multiple zeta value $\zeta(k_1,\cdots,k_r)$, denote by $N=k_1+\cdots+k_r$ and $r$ its weight  and depth respectively.  For $r=1$, these values correspond to the classical Riemann zeta values.   There are shuffle and stuffle products on the products of two multiple zeta values \cite{ikz}. 
  The theory of multiple zeta (zeta-star) values is related to modular forms \cite{GKZ}, mixed Tate motives \cite{brow, DG}, Grothendieck-Teichm\"{u}ller group \cite{furud}, knot theory \cite{BK}, fractal geometry \cite{lit,lir}, partial differential equations \cite{BGMR} and other subjects.
      
      Multiple zeta-star values are defined by
      \[
      \zeta^\star(k_1,\cdots,k_r)=\sum_{n_1\geq \cdots\geq n_r\geq 1}\frac{1}{n_1^{k_1}\cdots n_r^{k_r}}, k_1\geq 2,k_2,\cdots, k_{r}\geq 1.
      \]
    Denote by $\mathcal{Z}^{\star}$ the set of multiple zeta-star values. By the main results of Zhao \cite{zhao} and the relation between multiple zeta values and multiple zeta-star values, for complex numbers $$ \sum_{j=1}^i \mathrm{Re}(u_j)>1,i=1,\cdots,r-1,
    $$
    $$
      \sum_{j=1}^r \mathrm{Re}(u_j)>r , $$one can also define the multiple zeta-star function $$    \zeta^\star(u_1,u_2, \cdots,u_r) $$
       similarly.
      
      A fundamental problem in number theory is to determine the irrationality of many real numbers which arise naturally in mathematics and physics. For example, we know very little about the irrationality of Riemann zeta values, multiple zeta values and multiple zeta-star values. For more details, see \cite{ape},\cite{BR},\cite{beu},\cite{CDT}, \cite{fis},\cite{FSZ},\cite{LY}.
      
      Dirichlet's approximation theorem implies that: For $\alpha\in \mathbb{R}$, $\alpha\notin \mathbb{Q}$ if and only if
      \[
      \bigg{|}\,\alpha-\frac{p}{q}\,\bigg{|}<\frac{1}{q^2}
      \]
      for infinitely many integers $p,q$ and $(p,q)=1$.
      
      For rational approximation of real numbers, in 1924
       Khinchin \cite{khi} proved the following theorem, which is known as the Khinchin's zero-one law: \\
       Let $\psi:\mathbb{N}\rightarrow [0,+\infty)$ be a function such that $$q\psi(q)\geq (q+1)\psi(q+1), q\geq 1.$$
       Denote by $m$ the Lebesgue measure and $E_{\psi}$ the set of real numbers $\alpha\in [0,1]$ which satisfy 
       \[
         \bigg{|}\,\alpha-\frac{p}{q}\,\bigg{|}<\frac{\psi{(q)}}{q}       \]   
         for infinitely many $p,q$ and $(p,q)=1$. Then \\
         $(1)$ If $\sum_{q\geq 1}\psi(q)$ is convergent, then $m(E_{\psi})=0$.\\
         $(2)$ If $\sum_{q\geq 1}\psi(q)$ is divergent, then $m(E_{\psi})=1$.
         
         By the main result of Koukoulopoulos and Maynard \cite{KM}, the Khinchin's zero-one law  is generalized to the following theorem:   
           Let $\psi:\mathbb{N}\rightarrow [0,+\infty)$ be a function.        Denote by $m$ the Lebesgue measure and $A_{\psi}$ the set of real numbers $\alpha\in [0,1]$ which satisfy 
       \[
         \bigg{|}\,\alpha-\frac{p}{q}\,\bigg{|}<\frac{\psi{(q)}}{q}       \]   
         for infinitely many $p$ and $q$. Here $\varphi$ is the Euler's totient function.  \\
         $(1)$ If $\sum_{q\geq 1}\frac{\psi(q)\varphi(q)}{q}$ is convergent, then $m(A_{\psi})=0$.\\
         $(2)$ If $\sum_{q\geq 1}\frac{\psi(q)\varphi(q)}{q}$ is divergent, then $m(A_{\psi})=1$.    \\
         The statement $(2)$ in the above theorem is known as the Duffin-Schaeffer conjecture \cite{DS}. 
         Aistleitner, Borda and Hauke \cite{ABH} proved a quantitative version of the  Duffin-Schaeffer conjecture.          Recently, Fr$\ddot{\mathrm{u}}$hwirth and Hauke \cite{FH} established the Duffin-Schaeffer Conjecture for multiplicative
Diophantine approximation.         

           In \cite{lit}, the author found that there is a natural total order structure on the set of multiple zeta-star values (Theorem $1.2$ in \cite{lit}). Furthermore, the author showed that:  Define 
 \[
 \mathcal{T}=\Big{\{} (k_1,k_2, \cdots, k_r,\cdots)\,\Big{|}\,k_1\geq 2,k_i\geq 1, i\geq 2, k_s\geq 2 \mathrm{\;for \; some} \; s\geq 2 \;\mathrm{if} \;k_1=2    \Big{\}}.
 \]
 Then the map 
 \[
 \eta: \mathcal{T}\rightarrow (1, +\infty),
  \]
    \[
  {\bf k}=(k_1,k_2, \cdots, k_r,\cdots) \mapsto x=\mathop{\mathrm{lim}}_{r\rightarrow +\infty}\zeta^{\star}(k_1,k_2,\cdots,k_r)  \]
 is bijective.  We call this  map zeta-star correspondence.  The zeta-star correspondence should be compared with the theory of continued fractions, which gives a bijection between infinite sequences of positive integers and $(0,1)/\mathbb{Q}$.
  As a result, the set of multiple zeta-star values is a countable dense subset of $(1, +\infty)$.   In \cite{hmo}, M. Hirose, H. Murahara and T. Onozuka also found the bijective map $\eta$ independently.        
 
  The set of rational numbers is a countable dense subset of $\mathbb{R}$, while  the set of multiple zeta-star values is a countable dense subset of $(1, +\infty)$. In \cite{lit}, the author proved that
   \[
  \mathbb{Z}\; \bigcap\; \mathcal{Z}^\star=\emptyset.
  \]  An important open problem is to determine the following intersection set 
  \[
  \mathbb{Q}\; \bigcap\; \mathcal{Z}^\star.
  \]
  Thus there is a  natural question about the the set of multiple zeta-star values. For $\alpha>1$, can we give a criterion to determine $\alpha$ is a multiple zeta-star value or not?   
  In this paper, it is proved that
  \begin{Thm}\label{diri}
  For $\alpha>1$, then  $\alpha\notin \mathcal{Z}^\star$ if and only if  
  \[
   \bigg{|}\,\alpha-\zeta^\star(k_1,\cdots,k_r)\,\bigg{|}<\sum_{n_1\geq\cdots \geq n_r\geq n_{r+1}\geq 2} \frac{1}{n_1^{k_1}\cdots n_r^{k_r} n_{r+1}}    \]
   for infinitely many $r, k_1, \cdots,k_r$ with $r\rightarrow +\infty$.
  \end{Thm}
   Here we recall that 
   \[
   \sum_{n_1\geq\cdots \geq n_r\geq n_{r+1}\geq 2} \frac{1}{n_1^{k_1}\cdots n_r^{k_r} n_{r+1}}  = \zeta^\star(k_1,\cdots,k_r,1)-\zeta^\star(k_1,\cdots,k_r). \]  
   Thorem \ref{diri} can be viewed as the analogue of Dirichlet's approximation theorem  for the set of multiple zeta-star values. Roughly speaking, if a number is very close to infinitely many multiple zeta-star values, then this number is {\bf{not}} a multiple zeta-star value. Unlike the set of rational numbers, the set $\mathcal{Z}^\star$ is not an additive set.  For the proof of Theorem \ref{diri} we use the order structure of multiple zeta values and estimations of some  multiple series.
   
   Based on  Theorem \ref{diri} and the classical Diophantine approximation theory, there are the following questions: Are there any similar results like  Khinchin's zero-one law and Duffin-Schaeffer conjecture? We have the following result.
  
  \begin{Thm}\label{01}
  For a function $f: \mathbb{N}\rightarrow [1,+\infty)$. Denote by  $F_{f}$ the set of real numbers $\alpha\in (1,+\infty)$ which satisfy 
     \[
   \bigg{|}\,\alpha-\zeta^\star(k_1,\cdots,k_r)\,\bigg{|}<\sum_{n_1\geq\cdots \geq n_r\geq n_{r+1}\geq 2} \frac{1}{n_1^{k_1}\cdots n_r^{k_r} n_{r+1}^{f(r)}}    \]     for infinitely many $r, k_1, \cdots,k_r$ with $r\rightarrow +\infty$.
       \\
         $(i)$ If the series \[\sum_{r\geq 1}\frac{1}{2^{f(r)}}\] is convergent, then $m(F_{f})=0$.\\
         $(ii)$  If the series \[\sum_{r\geq 1}\frac{1}{2^{f(r)}}\] is divergent, then $m((1,+\infty)-F_{f})=0$.
  \end{Thm}
  
    Here we recall that 
   \[
   \sum_{n_1\geq\cdots \geq n_r\geq n_{r+1}\geq 2} \frac{1}{n_1^{k_1}\cdots n_r^{k_r} n_{r+1}^{f(r)}}  = \zeta^\star(k_1,\cdots,k_r,f(r))-\zeta^\star(k_1,\cdots,k_r). \]    Theorem \ref{01} is the zero-one law for the set of multiple zeta-star values. For the proof of $(i)$,  we will use the same analysis as the proof of Khinchin's zero-one law. While the proof of $(ii)$ is more tricky, we will use the order structure of multiple zeta-star values and the zeta-star correspondence, which is similar to the theory of continued fractions \cite{kh}.
  
    For a function $f: \mathbb{N}\rightarrow [1,+\infty)$. Denote by  $B_{f}$ the set of real numbers $\alpha\in [1,+\infty)$ which satisfy 
     \[
   \bigg{|}\,\alpha-\zeta^\star(k_1,\cdots,k_r)\,\bigg{|}<\frac{1}{2^{k_1+\cdots+k_r+f(r)}}    \]     for infinitely many $r, k_1, \cdots,k_r$ with $r\rightarrow +\infty$.
     If \[\sum_{r\geq 1}\frac{1}{2^{f(r)}}\] is convergent, as $B_f\subseteq F_f$, by Theorem \ref{01}, $(i)$, then $m(B_{f})=0$.  Now we propose the following conjecture.

  \begin{conj}\label{con}
 
           If \[\sum_{r\geq 1}\frac{1}{2^{f(r)}}\] is divergent,  
                    then $m((1,+\infty)-B_{f})=0$.  \end{conj}
                    
                    By the Theorem \ref{diri}, Theorem \ref{01} and Conjecture \ref{con}, one can see that the following kinds of  multiple series 
                    \[
                    \sum_{n_1\geq n_2\geq \cdots \geq n_r\geq 2}\frac{1}{n_1^{k_1}n_2^{k_2}\cdots n_r^{k_r}},                     \]
                    play an important role in the Diophantine approximation theory of multiple zeta-star values. In general cases, we have 
\begin{Thm}\label{m}
$(i)$ For $u>0$, 
\[
\sum_{\substack{  k_0+k_1+\cdots+k_r=K      \\  k_0,k_1,\cdots,k_r\geq 1      }} \sum_{n_1\geq \cdots n_r\geq 1}\frac{1}{(n_1+u-1)^{k_0} (n_1+u)^{k_1}(n_2+u)^{k_2}\cdots (n_r+u)^{k_r}}=\frac{1}{u^{K-r}};
\]
$(ii)$ For $m\geq 2$, one has 
\[
\sum_{\substack{  k_1+\cdots +k_r=K  \\k_1\geq 2,k_2,\cdots,k_r\geq 1      }} (k_1-1)  \sum_{n_1\geq n_2\geq \cdots \geq n_r\geq m}\frac{1}{n_1^{k_1}n_2^{k_2}\cdots n_r^{k_r}}<\frac{1}{(m-1)^{K-r}}.
\]
\end{Thm}
  
 Now  we give some remarks about the main results of this paper. 
  \begin{rem}
      The Grothdendieck period conjecture or the weight graded conjecture \cite{ikz} of multiple zeta values implies that every multiple zeta value is transcendental. Thus we  should have:
    \[
     \mathrm{If}\; \alpha\in\mathbb{Q}, \mathrm{then}\;\alpha\notin \mathcal{Z}^\star .
     \]
Can we use Theorem \ref{diri} to prove the above statement?  As
\[\mathbb{Q}=\bigcup_{m\geq 1}\frac{1}{m}\mathbb{Z},
\]
can we prove that: For some $m\geq 2$, if $\alpha \in \frac{1}{m}\mathbb{Z}$, then 
$\alpha\notin \mathcal{Z}^\star$? 
\end{rem}
  \begin{rem}
  By Theorem \ref{01}, if the series 
  \[
   \sum_{r\geq 1}\frac{1}{2^{f(r)}}
  \]
  is convergent, then $F_f$ is   a subset of $(1,+\infty)$ of Lebesgue measure $0$. How to calculate the Hausdorff dimension of the set $F_f$ in this case?
  \end{rem}
  \begin{rem} (Multiplicative Diophantine approximation)
  In this paper,we establish some classical Diophantine approximation type results for the set of multiple zeta-star values. Can we generalize these results to the multiplicative Diophantine approximation  of  multiple zeta-star values? For example, can we propose the Littlewood's conjecture in $k$ dimension \cite{CY} for the set of multiple zeta-star values?
  \end{rem}
  
  \begin{rem} (Additive Diophantine approximation)
  The sum of two rational numbers is also rational. But the sum of two multiple zeta-star values may not be a multiple zeta-star value.  As 
  $$\mathcal{Z}^\star+\mathcal{Z}^\star$$ is also a countable dense subset of the half line $(2,+\infty)$, 
  can we establish a criterion for $\alpha \in \mathcal{Z}^\star+\mathcal{Z}^\star$ or not? More generally, can we establish a criterion for $$\alpha \in  \underbrace{\mathcal{Z}^\star+ \mathcal{Z}^\star  +\cdots +\mathcal{Z}^\star}_{n}$$ or not?
  In fact, if one can show that 
  \[
  m\notin  \underbrace{\mathcal{Z}^\star+ \mathcal{Z}^\star  +\cdots +\mathcal{Z}^\star}_{n}
    \]
    for any integers  $m>n\geq 1$, then every multiple zeta-star values is irrational.
      \end{rem}

               \section{The order structure of multiple zeta-star values}\label{gea}
        In this section we review the order structure of multiple zeta-star values. For more details, see \cite{lit}.
Denote by 
\[
\mathcal{S}=\Big{\{}(k_1,\cdots,k_r)\,\Big{|}\,k_1\geq 2, k_2,\cdots,k_r\geq 1, r\geq1\Big{\}}.
\]
Define an order $\succ$ on $\mathcal{S}$ by
\[
(k_1,\cdots, k_{r},k_{r+1})\succ(k_1,\cdots, k_{r}), \forall\, (k_1,\cdots,k_r, k_{r+1})\in \mathcal{S}\]
and 
\[
(k_1,\cdots, k_{r})\succ(m_1,\cdots, m_{s})
\]
if $  k_i=m_i, 1\leq i\leq j, k_{j+1}<m_{j+1}$, for some $j\geq 0$.
 
 For $(k_1,\cdots,k_r)\in \mathcal{S}$, one has 
 \[
 \begin{split}
 &\;\;\;\; \zeta^\star(k_1,\cdots, k_{r},k_{r+1})\\
 &=\sum_{n_1\geq\cdots \geq n_r\geq n_{r+1}\geq 1}\frac{1}{n_1^{k_1}\cdots n_r^{k_r}n_{r+1}^{k_{r+1}}}\\
 &=\left(\sum_{\substack{n_1\geq\cdots \geq n_r\geq n_{r+1}\geq 1\\ n_{r+1}=1}}+\sum_{\substack{n_1\geq\cdots \geq n_r\geq n_{r+1}\geq 1\\ n_{r+1}\geq 2}}\right) \frac{1}{n_1^{k_1}\cdots n_r^{k_r}n_{r+1}^{k_{r+1}}}\\
 &=\zeta^\star(k_1,\cdots, k_{r})+\sum_{n_1\geq\cdots \geq n_r\geq n_{r+1}\geq 2}\frac{1}{n_1^{k_1}\cdots n_r^{k_r}n_{r+1}^{k_{r+1}}} 
 >\zeta^\star(k_1,\cdots, k_{r}). \end{split}
 \]
 
 In general cases, we have the following lemma.
 \begin{lem}\label{one}
 For $r\geq1$ and any fixed integer $m\geq 2$, one has \\
 $(i)$ \[
 \sum_{m\geq n_1\geq \cdots \geq n_r\geq 1}\frac{1}{n_1n_2\cdots n_r}< m;
 \]
 $(ii)$ \[
\mathop{\mathrm{lim}}_{r\rightarrow +\infty} \sum_{m\geq n_1\geq \cdots \geq n_r\geq 1}\frac{1}{n_1n_2\cdots n_r}= m;
 \]
 \end{lem}
\noindent{\bf Proof:}
$(i)$ For $r=1$, the statement $(i)$ is clear.  For $r>1$, one has 
\[
\frac{1}{n_1n_2\cdots n_r}=\frac{1}{n_1n_2\cdots n_{r-1}}, \;n_r=1;\]
\[
\frac{1}{n_1n_2\cdots n_r}<\frac{1}{n_1n_2\cdots n_{r-1}}, \;n_r>1.\]
Thus 
\[
\begin{split}
&\;\;\;\;\sum_{m\geq n_1\geq \cdots \geq n_r\geq 1}\frac{1}{n_1n_2\cdots n_r}<\sum_{m\geq n_1\geq \cdots \geq n_r\geq 1}\frac{1}{n_1n_2\cdots n_{r-1}}.\\
\end{split}
\]
One can check that
\[
\begin{split}
&\;\;\;\;\sum_{m\geq n_1\geq \cdots \geq n_r\geq 1}\frac{1}{n_1n_2\cdots n_{r-1}}\\
&=   \sum_{m\geq n_1\geq \cdots \geq n_{r-1}\geq 1}\frac{n_{r-1}}{n_1n_2\cdots n_{r-1}}     \\
&=   \sum_{m\geq n_1\geq \cdots \geq n_{r-1}\geq 1}\frac{1}{n_1n_2\cdots n_{r-2}}     \\
&\;\;\;\;\dots\;\;\;\;\dots\;\;\;\;\dots\\
&=\sum_{m\geq n_1\geq 1}1\\
&=m.
\end{split}
\]
So we have 
\[
\begin{split}
&\;\;\;\;\sum_{m\geq n_1\geq \cdots \geq n_r\geq 1}\frac{1}{n_1n_2\cdots n_r}<m.\\
\end{split}
\]
$(ii)$ From the following simple observation
\[
\begin{split}
&\;\;\;\;\sum_{m\geq n_1\geq \cdots \geq n_r\geq 1}\frac{1}{n_1n_2\cdots n_r}\\
&=\left( \sum_{\substack{m\geq n_1\geq \cdots \geq n_r\geq 1\\n_1=\cdots=n_r=1 }}+\sum_{\substack{m\geq n_1\geq \cdots \geq n_r\geq 1\\n_1\geq 2,   n_2=\cdots=n_r=1    }}+\cdots+ \sum_{\substack{m\geq n_1\geq \cdots \geq n_r\geq 1\\ n_r\geq 2      }} \right)\frac{1}{n_1n_2\cdots n_r}\\
&=1+\sum_{m\geq n_1\geq 2}\frac{1}{n_1}+ \cdots \sum_{m\geq n_1\geq\cdots\geq n_r\geq 2}\frac{1}{n_1n_2\cdots n_r}  ,    \\
\end{split}     
\]       
one has  
\[
\begin{split}
&\;\;\;\; \mathop{\mathrm{lim}}_{r\rightarrow +\infty} \sum_{m\geq n_1\geq \cdots \geq n_r\geq 1}\frac{1}{n_1n_2\cdots n_r}\\
&= 1+\sum_{m\geq n_1\geq 2}\frac{1}{n_1} +\sum_{m\geq n_1\geq n_2\geq 2}\frac{1}{n_1n_2} + \cdots \sum_{m\geq n_1\geq\cdots\geq n_r\geq 2}\frac{1}{n_1n_2\cdots n_r} +\cdots     \\
&=\prod_{n=2}^m \left( 1+\frac{1}{n}+\frac{1}{n^2}+\cdots+\frac{1}{n^r}+\cdots    \right)\\
&=m.
\end{split}
\]
$\hfill\Box$\\

If \[
(k_1,\cdots, k_{r})\succ(m_1,\cdots, m_{s})
\]
and  $  k_i=m_i, 1\leq i\leq j, k_{j+1}<m_{j+1}$, for some $j\geq 0$, then
\[
\begin{split}
&\;\;\;\;   \zeta^\star(m_1,\cdots, m_{s})
   \\
   &=\sum_{n_1\geq \cdots\geq n_j\geq n_{j+1}\geq \cdots\geq  n_s\geq 1}\frac{1}{n_1^{m_1}\cdots n_j^{m_j} n_{j+1}^{m_{j+1}}\cdots n_s^{m_s}}\\
   &\leq \sum_{n_1\geq \cdots\geq n_j\geq n_{j+1}\geq \cdots\geq  n_s\geq 1}\frac{1}{n_1^{m_1}\cdots n_j^{m_j} n_{j+1}^{m_{j+1}}n_{j+2}\cdots n_s}\\
   &<\sum_{n_1\geq \cdots\geq n_j\geq n_{j+1}\geq 1}\frac{1}{n_1^{m_1}\cdots n_j^{m_j} n_{j+1}^{m_{j+1}-1}}.\\
   \end{split}
   \]
    Here the last inequality follows from Lemma \ref{one}, $(i)$.
   Since $k_{j+1}<m_{j+1}$, 
   \[
   \sum_{n_1\geq \cdots\geq n_j\geq n_{j+1}\geq 1}\frac{1}{n_1^{m_1}\cdots n_j^{m_j} n_{j+1}^{m_{j+1}-1}}\leq \sum_{n_1\geq \cdots\geq n_j\geq n_{j+1}\geq 1}\frac{1}{n_1^{k_1}\cdots n_j^{k_j} n_{j+1}^{k_{j+1}}} \leq   \zeta^\star(k_1,\cdots, k_{r}).
 \]
  As a result,
  \[
  \zeta^\star(k_1,\cdots, k_{r})>  \zeta^\star(m_1,\cdots, m_{s}). \]
  
  Since $(\mathcal{S}, \succ)$ is a total order set, by induction we have:
  For $$(k_1,\cdots, k_r), (m_1,\cdots,m_s)\in\mathcal{S},$$
  \[
  (k_1,\cdots, k_r)\succ (m_1,\cdots,m_s)  \]
  if and only if 
  \[
    \zeta^\star(k_1,\cdots, k_r)> \zeta^\star(m_1,\cdots,m_s).  \]
    Furthermore, the order $\succ$ on $\mathcal{S}$ can be generalized to the order $\succ$ on $\mathcal{T}$. Similarly, one has: 
    For $${\bf{k}}=(k_1,\cdots, k_r,\cdots), {\bf{m}}=(m_1,\cdots,m_s,\cdots)\in\mathcal{T},$$
  \[
{\bf k}  \succ  {\bf m} \]
  if and only if 
  \[
    \eta({\bf k} )>\eta( {\bf m}).  \]    
    As an application of Lemma \ref{one}, $(ii)$, one can see that 
    \[
     \zeta^\star(k_1,\cdots, k_r)=\mathop{\mathrm{lim}}_{s\rightarrow +\infty}  \zeta^\star(k_1,\cdots, k_r+1, \{1\}^s) .  \]

   \section{Dirichlet approximation  of multiple zeta-star values}\label{sie}
   In this section, by using the order structure of multiple zeta-star values and zeta-star correspondence, we will prove Theorem \ref{diri}. 
   \begin{lem}\label{<} For $p\geq 2$, one has 
   \[
   \frac{3}{4}<\sum_{n_1\geq \cdots\geq n_p\geq 2}\frac{1}{n_1^2n_2\cdots n_p}<1.
   \]
   \end{lem}
   \noindent{\bf Proof:}
   As 
   \[
   \zeta^\star(2,\{1\}^{p-1})= \sum_{n_1\geq \cdots\geq n_p\geq 1}\frac{1}{n_1^2n_2\cdots n_p}=  \left(\sum_{\substack{n_1\geq \cdots\geq n_p\geq 1\\  n_p=1     }} + \sum_{\substack{n_1\geq \cdots\geq n_p\geq 1\\n_p\geq 2   }} \right) \frac{1}{n_1^2n_2\cdots n_p}, \]
   by the formula $\zeta^\star(2,\{1\}^{p-1})=p\zeta(p+1)$  (see \cite{ow},\cite{zlo}), we have 
   \[
   \begin{split}
   &\;\;\;\;   \sum_{n_1\geq \cdots\geq n_p\geq 2}\frac{1}{n_1^2n_2\cdots n_p}     \\
   &=  p\zeta(p+1)-(p-1)\zeta(p)      \\
   &=1-\sum_{n\geq 2} \frac{n(p-1)-p}{n^{p+1}}\\
   &<1
   \end{split}
   \] 
   for $p\geq 2$.
   Also, it is easy to check that 
   \[
   \frac{n(p-1)-p}{n^{p+1}}> \frac{np-(p+1)}{n^{p+2}}   \]
    for $n,p\geq 2$.
   Thus 
   \[
   \sum_{n_1\geq \cdots\geq n_p\geq 2}\frac{1}{n_1^2n_2\cdots n_p} \geq 2\zeta(3)-\zeta(2)>2\times 1.2-\frac{\pi^2}{6}>\frac{3}{4}.   \]
     
   $\hfill\Box$\\

    \begin{lem}\label{msv}
For  $k_1\geq 2,k_2,\cdots, k_r\geq 1$ and $n_1\geq \cdots\geq n_{r+1}\geq n_{r+2}\geq 2$, it follows that: \\
(i) If $n_1=2$, then
\[
\frac{1}{n_1^{k_1}\cdots n_r^{k_r}n_{r+1}n_{r+2}}= \frac{1}{(n_1+1)n_1n_2^{k_1-1}n_3^{k_2}\cdots n_{r+1}^{k_r}n_{r+2}}+\frac{1}{(n_1+2)(n_1+1)n_2n_3^{k_1-1}n_4^{k_2}\cdots n_{r+2}^{k_r}};
\]
(ii) If $n_1>2$, then
\[
\frac{1}{n_1^{k_1}\cdots n_r^{k_r}n_{r+1}n_{r+2}}< \frac{1}{(n_1+1)n_1n_2^{k_1-1}n_3^{k_2}\cdots n_{r+1}^{k_r}n_{r+2}}+\frac{1}{(n_1+2)(n_1+1)n_2n_3^{k_1-1}n_4^{k_2}\cdots n_{r+2}^{k_r}}.
\]
\end{lem}
 \noindent{\bf Proof:}
 If $n_1=2$, then $n_2=\cdots =n_{r+2}=2$. The identity in $(i)$ follows immediately  from
 \[
 \frac{1}{2^2}=\frac{1}{3\times 2}+\frac{1}{4\times 3}.
 \]
 If $n_1>2$, the inequality in $(ii)$ is equivalent to the following inequality
 \[
 \begin{split}
 &\;\;\;\;\frac{1}{n_1^{k_1}\cdots n_r^{k_r}n_{r+1}}\\
 &< \frac{1}{(n_1+1)n_1n_2^{k_1-1}n_3^{k_2}\cdots n_{r+1}^{k_r}}+\frac{1}{(n_1+2)(n_1+1)n_2n_3^{k_1-1}n_4^{k_2}\cdots n_{r+1}^{k_{r-1}}n_{r+2}^{k_r-1}}.\\
 \end{split} \tag{A}
 \]
 As $n_{r+2}\leq n_{r+1}$ and $k_r-1\geq 0$, it suffices to prove $(A)$ in case $n_{r+2}= n_{r+1}$. 
 
 For $n_{r+2}= n_{r+1}$, by multiplying $n_{r+1}$ on both sides, the inequality $(A)$ reduces to 
 \[
  \begin{split}
 &\;\;\;\;\frac{1}{n_1^{k_1}\cdots n_r^{k_r}}\\
 &< \frac{1}{(n_1+1)n_1n_2^{k_1-1}n_3^{k_2}\cdots n_r^{k_{r-1}} n_{r+1}^{k_r-1}}\\
 &\;\;\;\;+\frac{1}{(n_1+2)(n_1+1)n_2n_3^{k_1-1}n_4^{k_2}\cdots n_r^{k_{r-2}}n_{r+1}^{k_{r-1}+k_r-2}}.\\
 \end{split} \tag{B}
  \]
  Since $n_{r+1}\leq n_{r}$ and 
  \[
  k_r-1\geq 0,\]
  \[
   k_{r-1}+k_r-2\geq 0,
  \]
  it suffices to prove $(B)$ in case $n_{r+1}=n_{r}$. For $n_{r+1}=n_{r}$, the inequality $(B)$ reduces to 
   \[
  \begin{split}
 &\;\;\;\;\frac{1}{n_1^{k_1}\cdots n_{r}^{k_{r}}}\\
 &< \frac{1}{(n_1+1)n_1n_2^{k_1-1}n_3^{k_2}\cdots n_{r-1}^{k_{r-2}}n_r^{k_{r-1}+k_r-1} }\\
 &\;\;\;\;+\frac{1}{(n_1+2)(n_1+1)n_2n_3^{k_1-1}n_4^{k_2}\cdots n_{r-1}^{k_{r-3}}n_{r}^{k_{r-2}+k_{r-1}+k_r-2}}.\\
 \end{split} \tag{C}
  \]
  By multiplying $n_r^{k_r}$, the inequality $(C)$ reduces to 
   \[
  \begin{split}
 &\;\;\;\;\frac{1}{n_1^{k_1}\cdots n_{r-1}^{k_{r-1}}}\\
 &< \frac{1}{(n_1+1)n_1n_2^{k_1-1}n_3^{k_2}\cdots n_{r-1}^{k_{r-2}}n_r^{k_{r-1}-1} }+\frac{1}{(n_1+2)(n_1+1)n_2n_3^{k_1-1}n_4^{k_2}\cdots n_{r-1}^{k_{r-3}}n_{r}^{k_{r-2}+k_{r-1}-2}}.\\
 \end{split} \tag{D}
  \]
  Since $n_{r}\leq n_{r-1}$ and 
  \[
  k_{r-1}-1\geq 0,
  \]
  \[
  k_{r-2}+k_{r-1}-2\geq 0,
  \]
  it suffices to prove $(D)$ in case $n_r=n_{r-1}$.  For $n_{r}=n_{r-1}$, the inequality $(D)$ reduces to
   \[
  \begin{split}
 &\;\;\;\;\frac{1}{n_1^{k_1}\cdots n_{r-1}^{k_{r-1}}}\\
 &< \frac{1}{(n_1+1)n_1n_2^{k_1-1}n_3^{k_2}\cdots n_{r-2}^{k_{r-3}}n_{r-1}^{k_{r-2}+k_{r-1}-1} }\\
 &\;\;\;\;+\frac{1}{(n_1+2)(n_1+1)n_2n_3^{k_1-1}n_4^{k_2}\cdots n_{r-2}^{k_{r-4}}n_{r-1}^{k_{r-3}+k_{r-2}+k_{r-1}-2}}.\\
 \end{split} \tag{E}
  \]
  It is easy to see that the the inequalities $(C)$ and $(E)$ are of the same form.
  By repeating the above procedure for finite steps, it suffices to prove inequality $(A)$ in case $n_1=n_2=\cdots=n_{r+2}$.
          
  For \[
  n_1=n_2=\cdots=n_{r+2}>2  \]
   the inequality $(A)$ eventually reduces to
   \[
   \begin{split}
   &\;\;\;\;  \frac{1}{n_1^{k_1+\cdots+k_r+1} }         \\
   &<    \frac{1}{(n_1+1)n_1^{ k_1+\cdots+ k_r}}   +\frac{1}{ (n_1+2)(n_1+1)n_1^{k_1+\cdots+k_r-1}}    .  \\
   \end{split} \tag{F}
   \]
   By multiplying $n_1^{k_1+\cdots+k_r-1}$ on both sides, the inequality $(F)$ is equivalent to 
   \[
   \frac{1}{n_1^2}<\frac{1}{(n_1+1)n_1}+\frac{1}{(n_1+2)(n_1+1)}. \tag{G}
   \]
  As \[
  \frac{1}{(n_1+1)n_1}+\frac{1}{(n_1+2)(n_1+1)}=\frac{1}{n_1}-\frac{1}{n_1+2} =\frac{2}{n_1(n_1+2)}, \]
  the inequality $(G)$ is equivalent to 
  \[
  \frac{1}{n_1}<\frac{2}{n_1+2}. \tag{H}
  \]
  The inequality $(H)$ is equivalent to 
  \[
  n_1> \frac{n_1+2}{2}. \tag{I}
  \]
  By direct calculation, the inequality $(I)$ is equivalent to 
  \[n_1>2.\]
  As a result, the statement $(ii)$ is proved.
    $\hfill\Box$\\

 \begin{lem}\label{mzv}
For $k_1\geq 2,k_2,\cdots, k_r\geq 1$, we have 
\[
\begin{split}
&\sum_{n_1\geq \cdots\geq n_{r+2}\geq 2}\frac{1}{ n_1^{k_1}\cdots n_r^{k_r}n_{r+1} n_{r+2}       }<\sum_{n_1\geq \cdots\geq n_{r+1}\geq 2}  \frac{1}{ n_1^{k_1}\cdots n_r^{k_r}n_{r+1}     }  + \sum_{n_1\geq \cdots\geq n_{r}\geq 2}  \frac{1}{ n_1^{k_1}\cdots n_r^{k_r}   }  \\
\end{split}
\]
\end{lem}
 \noindent{\bf Proof:}
 It is clear that 
 \[
 \sum_{n\geq m}\frac{1}{(n+1)n}=\frac{1}{m}
 \]
 for $m\geq 1$. Thus
 \[
 \begin{split}
 &  \;\;\;\;     \sum_{n_1\geq \cdots\geq n_{r+1}\geq 2}  \frac{1}{ n_1^{k_1}\cdots n_r^{k_r}n_{r+1}     }     \\
 &=      \sum_{n_1\geq \cdots\geq n_{r+1}\geq 2}  \frac{1}{n_1}\cdot     \frac{1}{ n_1^{k_1-1}n_2^{k_2}\cdots n_r^{k_r}n_{r+1}     }       \\
 &=     \sum_{m_1\geq n_1\geq \cdots\geq n_{r+1}\geq 2}  \frac{1}{(m_1+1)m_1}\cdot     \frac{1}{ n_1^{k_1-1}n_2^{k_2}\cdots n_r^{k_r}n_{r+1}     }      \\
 &=  \sum_{n_1\geq n_2\geq \cdots\geq n_{r+2}\geq 2}      \frac{1}{ (n_1+1)n_1n_2^{k_1-1}n_3^{k_2}\cdots n_{r+1}^{k_r}n_{r+2}     }   \end{split}
 \]
 and 
 \[
 \begin{split}
 &  \;\;\;\;    \sum_{n_1\geq \cdots\geq n_{r}\geq 2}  \frac{1}{ n_1^{k_1}\cdots n_r^{k_r}   }     \\
 &=     \sum_{n_1\geq \cdots\geq n_{r}\geq 2}  \frac{1}{n_1}\cdot     \frac{1}{ n_1^{k_1-1}n_2^{k_2}\cdots n_r^{k_r}   }                 \\
 &=   \sum_{m_2\geq n_1\geq \cdots\geq n_{r+1}\geq 2}  \frac{1}{(m_2+1)m_2}\cdot     \frac{1}{ n_1^{k_1-1}n_2^{k_2}\cdots n_r^{k_r}    } \\
 &=\sum_{m_1\geq m_2\geq n_1\geq \cdots\geq n_{r}\geq 2}  \frac{1}{(m_1+2)(m_1+1)m_2}\cdot     \frac{1}{ n_1^{k_1-1}n_2^{k_2}\cdots n_r^{k_r}   }  \\
 &=\sum_{n_1\geq n_2\geq n_3\geq \cdots\geq n_{r+2}\geq 2}  \frac{1}{(n_1+2)(n_1+1)n_2}\cdot     \frac{1}{ n_3^{k_1-1}n_4^{k_2}\cdots n_{r+2}^{k_r}   }     \end{split}
 \] 
 By Lemma \ref{msv}, one has
 \[
 \begin{split}
 &\;\;\;\;  \sum_{n_1\geq \cdots\geq n_{r+2}\geq 2}\frac{1}{ n_1^{k_1}\cdots n_r^{k_r}n_{r+1} n_{r+2}       }     \\
 &<    \sum_{n_1\geq n_2\geq \cdots\geq n_{r+2}\geq 2}      \frac{1}{ (n_1+1)n_1n_2^{k_1-1}n_3^{k_2}\cdots n_{r+1}^{k_r}n_{r+2}     }\\
 &\;\;\;\;+   \sum_{n_1\geq n_2\geq n_3\geq \cdots\geq n_{r+2}\geq 2}       \frac{1}{ (n_1+2)(n_1+1)n_2 n_3^{k_1-1}n_4^{k_2}\cdots n_{r+2}^{k_r}   }  .     \\
 \end{split}
 \]
 So we have 
 \[
 \sum_{n_1\geq \cdots\geq n_{r+2}\geq 2}\frac{1}{ n_1^{k_1}\cdots n_r^{k_r}n_{r+1} n_{r+2}       }<\sum_{n_1\geq \cdots\geq n_{r+1}\geq 2}  \frac{1}{ n_1^{k_1}\cdots n_r^{k_r}n_{r+1}     }  + \sum_{n_1\geq \cdots\geq n_{r}\geq 2}  \frac{1}{ n_1^{k_1}\cdots n_r^{k_r}   } .
  \]
 $\hfill\Box$\\
   
       Now we are ready to prove Theorem \ref{diri}.

        For $\alpha>1$, if $\alpha\notin \mathcal{Z}^\star$, then by Theorem $1.3$ in \cite{lit} and 
   \[
     \zeta^\star(k_1,\cdots, k_r)=\mathop{\mathrm{lim}}_{s\rightarrow +\infty}  \zeta^\star(k_1,\cdots, k_r+1, \{1\}^s) ,  \]
     one has 
     \[
     \alpha=\eta((l_1,l_2,\cdots,l_r,\cdots))= \mathop{\mathrm{lim}}_{r\rightarrow +\infty} \zeta^\star(l_1,\cdots, l_r)     \]
     and $l_{r+1}\geq 2$ for infinitely many $r$. From the order structure of mutiple zeta-star values, we have 
     \[
       \zeta^\star(l_1,\cdots, l_r) < \alpha<   \zeta^\star(l_1,\cdots, l_r, l_{r+1}-1)  \leq \zeta^\star(l_1,\cdots, l_r, 1)   \]for infinitely many $r$.       
       As a result, one has
       \[
     \bigg{|}\alpha- \zeta^\star(l_1,\cdots, l_r) \bigg{|}<   \zeta^\star(l_1,\cdots, l_r, 1)-\zeta^\star(l_1,\cdots, l_r)= \sum_{n_1\geq \cdots \geq n_r\geq n_{r+1}\geq 2}\frac{1}{n_1^{l_1}\cdots n_r^{l_r} n_{r+1}}   \]for infinitely many $r,l_1,\cdots, l_r$ with $r\rightarrow +\infty$.  
     
     On the other hand, if 
      \[     \bigg{|}\alpha- \zeta^\star(k_1,\cdots, k_r) \bigg{|}< \sum_{n_1\geq \cdots \geq n_r\geq n_{r+1}\geq 2}  \frac{1}{n_1^{k_1}\cdots n_r^{k_r} n_{r+1}} \tag{1}  \]     for infinitely many $r,k_1,\cdots, k_r$   with $r\rightarrow +\infty$, then 
     \[
     \alpha<  \zeta^\star(k_1,\cdots, k_r,1 )      \]
     for infinitely many $r,k_1,\cdots, k_r$  with $r\rightarrow +\infty$.   
       
 Assuming that $\alpha\in \mathcal{Z}^\star$,
then $\alpha=\zeta^\star(l_1,l_2,\cdots,l_s)$ for some $(l_1,l_2,\cdots, l_s)\in \mathcal{S}$. 
By the order structure of multiple zeta-star values, one has 
\[
(k_1,k_2,\cdots, k_r,1)\succ (l_1,l_2,\cdots, l_s)
\]
for infinitely many $r,k_1,\cdots, k_r$  with $r\rightarrow +\infty$. There are the following two cases:\\
$(i)$ \[
(k_1,k_2,\cdots, k_r,1)\succ (k_1,k_2,\cdots, k_r)\succ (l_1,l_2,\cdots, l_s)
\]
for infinitely many $r,k_1,\cdots, k_r$  with $r\rightarrow +\infty$;\\
$(ii)$    \[
(k_1,k_2,\cdots, k_r,1)\succ (l_1,l_2,\cdots, l_s)\succ (k_1,k_2,\cdots, k_r)\]
for infinitely many $r,k_1,\cdots, k_r$ with $r\rightarrow +\infty$.
If $(1)$ and $(ii)$ is true, by the order structure of $\mathcal{S}$, one has $r<s$ and $k_i=l_i, 1\leq i\leq r$
for infinitely many $r,k_1,\cdots, k_r$. As $(l_1,l_2, \cdots,l_s)$ is fixed, it is not possible.

If $(1)$ and $(i)$ is true, for $r>s$ there are the following two cases:\\
$(iii)$ \[
(k_1,k_2,\cdots, k_r)\succ (k_1,k_2,\cdots, k_s)\succ (l_1,l_2,\cdots, l_s)
\]
for infinitely many $r,k_1,\cdots, k_r$  with $r\rightarrow +\infty$;\\
$(iv)$    \[
(k_1,k_2,\cdots, k_r)\succ (k_1,k_2,\cdots, k_s)=(l_1,l_2,\cdots, l_s)\]
for infinitely many $r,k_1,\cdots, k_r$  with $r\rightarrow +\infty$. For $r>s$, if $(1)$, $(i)$  is true, then 
\[
\begin{split}
&\;\;\;\;\bigg{|}  \alpha- \zeta^\star(k_1,\cdots, k_r)   \bigg{|}\\
&\geq  \zeta^\star(k_1,\cdots, k_r)- \zeta^\star(k_1,\cdots, k_s)\\
 &\geq \sum_{n_1\geq \cdots \geq n_{s+1}\geq 2}  \frac{1}{n_1^{k_1}\cdots  n_{s+1}^{k_{s+1}}}+\sum_{n_1\geq \cdots \geq n_{s+2}\geq 2}  \frac{1}{n_1^{k_1}\cdots  n_{s+2}^{k_{s+2}}} +\cdots + \sum_{n_1\geq \cdots \geq n_{r}\geq 2}  \frac{1}{n_1^{k_1}\cdots  n_{r}^{k_{r}}},
\end{split}
\]
\[
\bigg{|}  \alpha- \zeta^\star(k_1,\cdots, k_r)   \bigg{|}<\sum_{n_1\geq \cdots \geq n_r\geq n_{r+1}\geq 2}  \frac{1}{n_1^{k_1}\cdots n_r^{k_r} n_{r+1}} \]
for infinitely many $r,k_1,\cdots, k_r$  with $r\rightarrow +\infty$. 
As a result, 
\[
\begin{split}
&\;\;\;\;   \sum_{n_1\geq \cdots \geq n_{s+1}\geq 2}  \frac{1}{n_1^{k_1}\cdots  n_{s+1}^{k_{s+1}}}+\sum_{n_1\geq \cdots \geq n_{s+2}\geq 2}  \frac{1}{n_1^{k_1}\cdots  n_{s+2}^{k_{s+2}}} +\cdots + \sum_{n_1\geq \cdots \geq n_{r}\geq 2}  \frac{1}{n_1^{k_1}\cdots  n_{r}^{k_{r}}}         \\
&<\sum_{n_1\geq \cdots \geq n_r\geq n_{r+1}\geq 2}  \frac{1}{n_1^{k_1}\cdots n_r^{k_r} n_{r+1}}\end{split}
\]
for infinitely many $r,k_1,\cdots, k_r$  with $r>s$ and  $r\rightarrow +\infty$. 
In the above statement, if $k_j\geq 2$ for some $r\geq j\geq s+2$, then
\[
\begin{split}
&\;\;\;\;  \sum_{n_1\geq \cdots \geq n_r\geq n_{r+1}\geq 2}  \frac{1}{n_1^{k_1}\cdots n_r^{k_r} n_{r+1}}       \\
&= \sum_{n_1\geq \cdots \geq n_r\geq n_{j-1}\geq 2} \left( \frac{1}{n_1^{k_1}\cdots n_{j-1}^{k_{j-1}}} \cdot \sum_{n_{j-1}\geq n_j\geq \cdots\geq n_{r+1}\geq 2        }  \frac{1}{n_j^{k_j}\cdots n_r^{k_r} n_{r+1}}    \right)     \\
&\leq  \sum_{n_1\geq \cdots \geq n_r\geq n_{j-1}\geq 2} \left( \frac{1}{n_1^{k_1}\cdots n_{j-1}^{k_{j-1}}} \cdot \sum_{n_{j-1}\geq n_j\geq \cdots\geq n_{r+1}\geq 2        }  \frac{1}{n_j^{2}n_{j+1}\cdots n_r n_{r+1}}    \right)     \\
&<\sum_{n_1\geq \cdots \geq n_r\geq n_{j-1}\geq 2} \left( \frac{1}{n_1^{k_1}\cdots n_{j-1}^{k_{j-1}}} \cdot \sum_{+\infty>n_j\geq \cdots\geq n_{r+1}\geq 2        }  \frac{1}{n_j^{2}n_{j+1}\cdots n_r n_{r+1}}    \right)     \\
&<\sum_{n_1\geq \cdots \geq n_r\geq n_{j-1}\geq 2}  \frac{1}{n_1^{k_1}\cdots n_{j-1}^{k_{j-1}}} .\end{split}
\]
Here the last inequality follows from Lemma \ref{<}. The above inequality contradicts the previous statement. So we have 
$k_j=1$ for $r\geq j\geq s+2$ and 
\[
\begin{split}
&\;\;\;\;   \sum_{n_1\geq \cdots \geq n_{s+1}\geq 2}  \frac{1}{n_1^{k_1}\cdots  n_{s+1}^{k_{s+1}}}+\sum_{n_1\geq \cdots \geq n_{s+2}\geq 2}  \frac{1}{n_1^{k_1}\cdots n_{s+1}^{k_{s+1}} n_{s+2}} +\cdots \\
&+ \sum_{n_1\geq \cdots \geq n_{r}\geq 2}  \frac{1}{n_1^{k_1}\cdots n_{s+1}^{k_{s+1}} n_{s+2}\cdots  n_{r}}         
<\sum_{n_1\geq \cdots \geq n_r\geq n_{r+1}\geq 2}  \frac{1}{n_1^{k_1}\cdots n_{s+1}^{k_{s+1}} n_{s+2}\cdots  n_{r} n_{r+1}}\end{split} \tag{2}
\]
for infinitely many $r,k_1,\cdots, k_{s+1}$  with $r>s$ and  $r\rightarrow +\infty$.  By Lemma \ref{mzv}, one has 
\[
\begin{split}
&\;\;\;\; \sum_{n_1\geq \cdots \geq n_r\geq n_{r+1}\geq 2}  \frac{1}{n_1^{k_1}\cdots n_{s+1}^{k_{s+1}} n_{s+2}\cdots  n_{r} n_{r+1}}\\
&<\sum_{n_1\geq \cdots \geq n_{r-1}\geq n_{r}\geq 2}  \frac{1}{n_1^{k_1}\cdots n_{s+1}^{k_{s+1}} n_{s+2}\cdots  n_{r} }+ \sum_{n_1\geq \cdots \geq n_{r-2}\geq n_{r-1}\geq 2}  \frac{1}{n_1^{k_1}\cdots n_{s+1}^{k_{s+1}} n_{s+2}\cdots  n_{r-1} }. \\
\end{split}
\]
The above inequality contradicts the inequalities $(2)$. So we have $\alpha\notin \mathcal{Z}^\star$.

 As a result, Theorem \ref{diri} is proved.

\section{Zero-one law of multiple zeta-star values      }
     For a function $$f: \mathbb{N}\rightarrow [1,+\infty)$$ and $$k_1\geq 2, k_2,\cdots, k_r\geq 1.$$Define $\mathcal{U}_{k_1,\cdots,k_r}(f)$ and $\mathcal{V}_r(f)$ as 
     \[
     \mathcal{U}_{k_1,\cdots,k_r}(f)=\bigg{\{} x\, \bigg{|}\,x>1, \Big{|}x-\zeta^\star(k_1,\cdots,k_r)\Big{|}<\sum_{n_1\geq\cdots \geq n_r\geq n_{r+1}\geq 2} \frac{1}{n_1^{k_1}\cdots n_r^{k_r} n_{r+1}^{f(r)}}      \bigg{\}}   ,    \]
     \[
     \mathcal{V}_r(f)=\bigcup_{k_1\geq 2,k_2,\cdots, k_r\geq 1} \mathcal{U}_{k_1,\cdots,k_r}(f) .        \]
   By definition,  $F_{f}$ is the set of real numbers $\alpha\in [1,+\infty)$ which satisfy 
     \[
   \bigg{|}\,\alpha-\zeta^\star(k_1,\cdots,k_r)\,\bigg{|}<\sum_{n_1\geq\cdots \geq n_r\geq n_{r+1}\geq 2} \frac{1}{n_1^{k_1}\cdots n_r^{k_r} n_{r+1}^{f(r)}}    \]     for infinitely many $r, k_1, \cdots,k_r$ with $r\rightarrow +\infty$.     Thus one has
   \[
   F_f=\bigcap_{s\geq 1} \bigcup_{r\geq s} \mathcal{V}_r(f).   \]

   As a result, Theorem \ref{01} is equivalent to the following theorem.
   \begin{Thm}\label{eq01}
     $(i)$ If \[\sum_{r\geq 1}\frac{1}{2^{f(r)}}\] is convergent, then \[m\left( \bigcap_{s\geq 1} \bigcup_{r\geq s} \mathcal{V}_r(f) \right)=0.\]
         $(ii)$  If \[\sum_{r\geq 1}\frac{1}{2^{f(r)}}\] is divergent, then \[m\left((1,+\infty)-\bigcap_{s\geq 1} \bigcup_{r\geq s} \mathcal{V}_r(f)  \right)=0.\]
   \end{Thm}
   
    The following lemma will be used in the proof of Theorem \ref{eq01}, $(i)$.
        \begin{lem}\label{mone}
    For $r\geq t+2$ and $u\geq 1$, one has 
    \[
    \sum_{n_1\geq \cdots\geq n_{r+1}\geq 2}\frac{1}{n_1(n_1-1)(n_2-1)\cdots (n_r-1)}\cdot \frac{1}{n_t n_{r+1}^u}<\frac{\zeta^\star\left( 2,\{1\}^{t-1}    \right)}{2^u}.
    \]
    \end{lem}
      \noindent{\bf Proof:}
      By multiplying $2^u$ on both sides, the above inequality is equivalent to 
      \[
        \sum_{n_1\geq \cdots\geq n_{r+1}\geq 2}\frac{1}{n_1(n_1-1)(n_2-1)\cdots (n_r-1)}\cdot \frac{1}{n_t }\left( \frac{2}{n_{r+1}} \right)^u<\zeta^\star\left( 2,\{1\}^{t-1}    \right).      \]
        It is clear that
        \[
        \begin{split}
        &\;\;\;\;   \sum_{n_1\geq \cdots\geq n_{r+1}\geq 2}\frac{1}{n_1(n_1-1)(n_2-1)\cdots (n_r-1)}\cdot \frac{1}{n_t }\left( \frac{2}{n_{r+1}} \right)^u     \\
        &<    \sum_{n_1\geq \cdots\geq n_{r+1}\geq 2}\frac{1}{n_1(n_1-1)(n_2-1)\cdots (n_r-1)}\cdot \frac{1}{n_t }        \\
        &\leq      \sum_{n_1\geq \cdots\geq n_{t}\geq 2}\frac{1}{n_1(n_1-1)(n_2-1)\cdots (n_t-1)}\cdot \frac{1}{n_t } \left( \sum_{n_t\geq n_{t+1}\geq \cdots\geq n_{r+1}\geq 2}   \frac{1}{(n_{t+1}-1)\cdots (n_r-1)} \right)   .             \\
        \end{split}
        \]
        As 
        \[
        \begin{split}
        &\;\;\;\;    \sum_{n_t\geq n_{t+1}\geq \cdots\geq n_{r+1}\geq 2}   \frac{1}{(n_{t+1}-1)\cdots (n_r-1)}    \\
        &=   \sum_{n_t\geq n_{t+1}\geq \cdots\geq n_{r}\geq 2}   \frac{n_r-1}{(n_{t+1}-1)\cdots (n_r-1)}       \\
        &=  \sum_{n_t\geq n_{t+1}\geq \cdots\geq n_{r}\geq 2}   \frac{1}{(n_{t+1}-1)\cdots (n_{r-1}-1)}     \\
        &\cdots\\
        &=\sum_{n_t\geq n_{t+1}\geq 2} 1\\
        &= n_t-1,
        \end{split}
        \]
        one has 
           \[
        \begin{split}
        &\;\;\;\;   \sum_{n_1\geq \cdots\geq n_{r+1}\geq 2}\frac{1}{n_1(n_1-1)(n_2-1)\cdots (n_r-1)}\cdot \frac{1}{n_t }\left( \frac{2}{n_{r+1}} \right)^u     \\
        &\leq      \sum_{n_1\geq \cdots\geq n_{t}\geq 2}\frac{1}{n_1(n_1-1)(n_2-1)\cdots (n_t-1)}\cdot \frac{n_t-1}{n_t }            \\
        &<    \sum_{n_1\geq \cdots\geq n_{t}\geq 2}\frac{1}{n_1(n_1-1)(n_2-1)\cdots (n_t-1)}          \\
          &<    \sum_{n_1\geq \cdots\geq n_{t}\geq 2}\frac{1}{(n_1-1)^2(n_2-1)\cdots (n_t-1)} =\zeta^\star\left(2,\{1\}^{t-1}\right).         \\
        \end{split}
        \]
        In a word, the lemma is proved.      $\hfill\Box$\\
  \noindent{\bf Proof of Theorem \ref{eq01} $(i)$:}
For any $s\geq 1$, it is clear that
  \[
  m\left( \bigcap_{s\geq 1} \bigcup_{r\geq s} \mathcal{V}_r(f) \right)\leq m\left(  \bigcup_{r\geq s} \mathcal{V}_r(f) \right) .   \]
  Thus it suffices to show that 
  \[
     \lim_{s\rightarrow +\infty } m\left(  \bigcup_{r\geq s} \mathcal{V}_r(f) \right)=0.  \tag{1} \]
     Beware that $  \bigcup\limits_{r\geq s} \mathcal{V}_r(f)  $ is an unbounded subset of $(1,+\infty)$. The formula $(1)$ is equivalent to 
      \[
     \lim_{s\rightarrow +\infty } m\left(  \bigcup_{r\geq s} \mathcal{V}_r(f) \bigcap \,(1, \zeta^\star(2,\{1\}^{t-1})\right)=0. \tag{2} \] for 
     any $t\geq 1$.
         For $s\geq t+3$ We have 
     \[
     \begin{split}
     & \;\;\;\;   m\left(  \bigcup_{r\geq s} \mathcal{V}_r(f) \bigcap \,(1, \zeta^\star(2,\{1\}^{t-1}) \right)     \\
     &\leq \sum_{r\geq s} m\left(  \mathcal{V}_r(f) \bigcap \,(1, \zeta^\star(2,\{1\}^{t-1}) \right)    \\
     &\leq  \sum_{r\geq s}  \sum_{k_1\geq 2, k_2,\cdots,k_r\geq 1}  m\left(\mathcal{U}_{k_1,\cdots,k_r}(f) \bigcap \,(1, \zeta^\star(2,\{1\}^{t-1}) \right) . \\
     \end{split}
     \]
     By definition, one has 
     \[
     \zeta^\star(k_1,\cdots,k_r)+\sum_{n_1\geq \cdots\geq n_r\geq n_{r+1}\geq 2} \frac{1}{n_1^{k_1}\cdots n_r^{k_r}n_{r+1}^{f(r)}}= \zeta^\star(k_1,\cdots,k_r,f(r) )     \]
     and 
     \[
     \begin{split}
     &\;\;\;\; \mathcal{U}_{k_1,\cdots,k_r}(f)           \\
     &= \left(   \zeta^\star(k_1,\cdots,k_r)-\sum_{n_1\geq \cdots\geq n_r\geq n_{r+1}\geq 2} \frac{1}{n_1^{k_1}\cdots n_r^{k_r}n_{r+1}^{f(r)}}      ,\;\zeta^\star(k_1,\cdots,k_r,f(r) )  \right).
     \end{split}
     \]
     As $f(r)\geq 1$ and $r\geq s\geq t+3$, there are two cases. \\ 
     $(i)$  If $k_r\geq 2$, then
     \[
     \begin{split}
     &\;\;\;\; \sum_{n_1\geq \cdots\geq n_r\geq n_{r+1}\geq 2} \frac{1}{n_1^{k_1}\cdots n_r^{k_r}n_{r+1}^{f(r)}}  \\
     &=   \sum_{n_1\geq \cdots\geq n_{r-1}\geq 2} \frac{1}{n_1^{k_1}\cdots n_{r-1}^{k_{r-1}}} \cdot \left( \sum_{n_{r-1}\geq n_r\geq n_{r+1}\geq 2}     \frac{1}{n_r^{k_r}n_{r+1}^{f(r)}}    \right) \\
     &<    \sum_{n_1\geq \cdots\geq n_{r-1}\geq 2} \frac{1}{n_1^{k_1}\cdots n_{r-1}^{k_{r-1}}} \cdot \left( \sum_{+\infty >n_r\geq n_{r+1}\geq 2}     \frac{1}{n_r^{k_r}n_{r+1}^{f(r)}}    \right)      \\     
        &\leq     \sum_{n_1\geq \cdots\geq n_{r-1}\geq 2} \frac{1}{n_1^{k_1}\cdots n_{r-1}^{k_{r-1}}} \cdot \left( \sum_{+\infty >n_r\geq n_{r+1}\geq 2}     \frac{1}{n_r^{2}n_{r+1}}    \right)      \\   
              &\leq     \sum_{n_1\geq \cdots\geq n_{r-1}\geq 2} \frac{1}{n_1^{k_1}\cdots n_{r-1}^{k_{r-1}}}  .  
     \end{split}
     \]
     Here the last inequality follows from Lemma \ref{<}. \\
     $(ii)$ If $k_r=1$,   by Lemma \ref{mzv},
            \[
     \begin{split}
     &\;\;\;\; \sum_{n_1\geq \cdots\geq n_r\geq n_{r+1}\geq 2} \frac{1}{n_1^{k_1}\cdots n_r^{k_r}n_{r+1}^{f(r)}}  \\
     &\leq \sum_{n_1\geq \cdots\geq n_r\geq n_{r+1}\geq 2} \frac{1}{n_1^{k_1}\cdots n_{r-1}^{k_{r-1}} n_rn_{r+1}}  \\
     &<   \sum_{n_1\geq \cdots\geq n_r\geq 2} \frac{1}{n_1^{k_1}\cdots n_{r-1}^{k_{r-1}} n_r}   +    \sum_{n_1\geq \cdots\geq n_{r-1}\geq 2} \frac{1}{n_1^{k_1}\cdots n_{r-1}^{k_{r-1}} }    \\
     &<   \sum_{n_1\geq \cdots\geq n_r\geq 2} \frac{1}{n_1^{k_1}\cdots n_{r-1}^{k_{r-1}} n_r^{k_r}}   +    \sum_{n_1\geq \cdots\geq n_{r-1}\geq 2} \frac{1}{n_1^{k_1}\cdots n_{r-1}^{k_{r-1}} }    \\
     \end{split}
     \] 
     From the above calculations, it is clear that 
       \[
     \begin{split}
     &\;\;\;\; \mathcal{U}_{k_1,\cdots,k_r}(f)           \subseteq \left(   \zeta^\star(k_1,\cdots,k_{r-2}) ,\;\zeta^\star(k_1,\cdots,k_r,f(r) )  \right).
     \end{split}
     \]
     For $r\geq s\geq t+3$, if 
     \[
     (k_1,\cdots, k_{r-2}) \succ (2, \{1\}^{t-1}),
        \]
        then by the order structure of multiple zeta-star values,
        \[
        \mathcal{U}_{k_1,\cdots,k_r}(f)  \bigcap \left( 1, \zeta^\star(2,\{1\}^{t-1}   \right) = \emptyset.       \]
        As 
        \[
         (k_1,\cdots, k_{r-2}) \succ (2, \{1\}^{t-1})  \Leftrightarrow k_1=2,k_2=\cdots, k_t= 1,       \]
         so   \[
     \begin{split}
     & \;\;\;\;   m\left(  \bigcup_{r\geq s} \mathcal{V}_r(f) \bigcap \,(1, \zeta^\star(2,\{1\}^{t-1}) \right)     \\
          &\leq  \sum_{r\geq s}  \sum_{k_1\geq 2, k_2,\cdots,k_r\geq 1}  m\left(\mathcal{U}_{k_1,\cdots,k_r}(f) \bigcap \,(1, \zeta^\star(2,\{1\}^{t-1}) \right) \\
          &\leq    \sum_{r\geq s}  \sum_{\substack{k_1\geq 2, k_2,\cdots,k_r\geq 1\\(k_1,\cdots,k_t)\neq (2, \{1\}^{t-1})   }}  m\left(\mathcal{U}_{k_1,\cdots,k_r}(f)  \right)     \\
     \end{split}
     \]
     In conclusion, for $s\geq t+3$, one has 
     \[
     \begin{split}
      & \;\;\;\;   m\left(  \bigcup_{r\geq s} \mathcal{V}_r(f) \bigcap \,(1, \zeta^\star(2,\{1\}^{t-1}) \right)     \\
          &\leq    \sum_{r\geq s}  \sum_{\substack{k_1\geq 2, k_2,\cdots,k_r\geq 1\\(k_1,\cdots,k_t)\neq (2, \{1\}^{t-1})   }}  m\left(\mathcal{U}_{k_1,\cdots,k_r}(f)  \right)     \\
     &\leq \sum_{r\geq s}  \sum_{\substack{k_1\geq 2, k_2,\cdots,k_r\geq 1\\   (k_1,\cdots,k_t)\neq (2, \{1\}^{t-1})      }}  \sum_{n_1\geq\cdots \geq n_r\geq n_{r+1}\geq 2} \frac{2}{n_1^{k_1}\cdots n_r^{k_r} n_{r+1}^{f(r)}}        \\
     &< \sum_{r\geq s}\left(  \sum_{{ k_1}\geq 3, k_2,\cdots,k_r\geq 1} +  \sum_{k_1\geq 2, k_2\geq 2,k_3,\cdots,k_r\geq 1}+\cdots+ \sum_{k_1\geq 2, k_2, \cdots,k_{t-1}\geq 1, k_t\geq 2,k_{t+1},\cdots,k_r\geq 1}\right)\\
    &\;\;\;\;\sum_{n_1\geq\cdots \geq n_r\geq n_{r+1}\geq 2} \frac{2}{n_1^{k_1}\cdots n_r^{k_r} n_{r+1}^{f(r)}}        \\
    &\leq  \sum_{r\geq s}   \sum_{n_1\geq\cdots \geq n_r\geq n_{r+1}\geq 2} 
    \Bigg{(}  \sum_{{ k_1}\geq 3, k_2,\cdots,k_r\geq 1} +  \sum_{k_1\geq 2, k_2\geq 2,k_3,\cdots,k_r\geq 1}+ \cdots\\
    &\;\;\;\;+\sum_{k_1\geq 2, k_2, \cdots,k_{t-1}\geq 1, k_t\geq 2,k_{t+1},\cdots,k_r\geq 1}\Bigg{)}   \frac{2}{n_1^{k_1}\cdots n_r^{k_r} n_{r+1}^{f(r)}}      \\
    &\leq   \sum_{r\geq s}   \sum_{n_1\geq\cdots \geq n_{r+1}\geq 2} \Bigg{(}   \frac{1}{n_1^2(n_1-1)(n_2-1)\cdots (n_r-1)}+ \frac{1}{n_1(n_1-1)n_2(n_2-1)(n_3-1)\cdots (n_r-1)}\\
    &\;\;\;\;+\cdots+ \frac{1}{   n_1(n_1-1)(n_2-1)\cdots (n_{t-1}-1)n_t(n_t-1)(n_{t+1}-1)\cdots (n_r-1)  }\Bigg{)}\cdot \frac{2}{n_{r+1}^{f(r)}}\\
    &\leq   \sum_{r\geq s}   \sum_{n_1\geq\cdots \geq n_{r+1}\geq 2}  \frac{1}{n_1(n_1-1)(n_2-1)\cdots (n_r-1)}  \Bigg{(} \frac{1}{n_1}+ \frac{1}{n_2}+\cdots+ \frac{1}{   n_t  }\Bigg{)}\cdot \frac{2}{n_{r+1}^{f(r)}}\\
     &\leq   \sum_{r\geq s}   \sum_{n_1\geq\cdots \geq n_{r+1}\geq 2}  \frac{1}{n_1(n_1-1)(n_2-1)\cdots (n_r-1)}  \cdot \frac{t}{   n_t  }\cdot \frac{2}{n_{r+1}^{f(r)}}\\
      \end{split}
     \]
     Here the last inequality follows from the fact that: If $n_1\geq \cdots \geq n_t\geq 2$, 
     \[
     \frac{1}{n_1}+ \frac{1}{n_2}+\cdots+ \frac{1}{   n_t  } \leq \frac{t}{   n_t  }   .       \]
     By Lemma \ref{mone}, for $s\geq t+3$, we have 
     \[
     \begin{split}
     &\;\;\;\;  m\left(  \bigcup_{r\geq s} \mathcal{V}_r(f) \bigcap \,(1, \zeta^\star(2,\{1\}^{t-1}) \right)         \\
     &< \sum_{r\geq s}   \frac{2t\cdot \zeta^\star\left(2,\{1\}^{t-1}\right)}{2^{f(r)}}  = 2t\cdot \zeta^\star\left(2,\{1\}^{t-1}\right)\cdot \sum_{r\geq s} \frac{1}{2^{f(r)}}. \\
     \end{split}
     \]
     Since the series 
     \[
     \sum_{r\geq 1}\frac{1}{2^{f(r)}}
     \]
     is convergent, one has 
           \[
     \lim_{s\rightarrow +\infty } m\left(  \bigcup_{r\geq s} \mathcal{V}_r(f) \bigcap \,(1, \zeta^\star(2,\{1\}^{t-1})\right)=0 \] for any $t\geq 1$. Thus 
             \[
     \lim_{s\rightarrow +\infty } m\left(  \bigcup_{r\geq s} \mathcal{V}_r(f)\right)=0 .    \]
     So we have 
    $
     m\left( \bigcap\limits_{s\geq 1} \bigcup\limits_{r\geq s} \mathcal{V}_r(f) \right)=0.
$  $\hfill\Box$\\
      
      The following lemma is a standard exercise in analysis.       \begin{lem}\label{prod}
    For  $f:\mathbb{N}\rightarrow [1,+\infty)$, if the series
 \[
 \sum_{r\geq 1}\frac{1}{2^{f(r)}}
 \]
 is divergent, then  
 \[
 \lim_{l\rightarrow +\infty}\prod_{r=1}^l \left(  1-\frac{1}{2^{f(r)}}  \right)=0.
 \]    \end{lem}
 \noindent{\bf Proof:}
The above statement follows immediately from the following inequality
 \[
 -\frac{x}{2}>\mathrm{log}\;(1-x)>-x,\; x\in (0,1).
 \]
 $\hfill\Box$\\  
 
 For $k_1\geq 2,k_2,\cdots,k_{r}\geq 1$ and 
 \[
 (k_1,k_2,\cdots,k_r)\neq (2,\{1\}^{r-1}),
 \]
 define ${Z}_{k_1,\cdots,k_r}$ as follows:\\
 \[
Z_{k_1,\cdots,k_r}=\left(  \zeta^\star(k_1,\cdots,k_{r-1},k_r), \zeta^\star(k_1,\cdots,k_{r-1},k_r-1)   \right]
  \]
  if $k_r\geq 2$;
  \[
  Z_{k_1,\cdots,k_r}=\left(  \zeta^\star(k_1,\cdots,k_{r-1},k_r), \zeta^\star(k_1,\cdots,k_{i-1},k_i-1)   \right]  \]
  if $$k_i\geq 2, k_{i+1}=\cdots=k_r=1.$$
  It is clear that 
  \[
Z_{k_1,\cdots,k_r}\bigcap Z_{l_1,\cdots,l_r} =\emptyset   \]
  for $(k_1,\cdots,k_r)\neq(l_1,\cdots,l_r)$.
   As 
  \[
  \zeta^\star(k_1,\cdots,k_r)=\lim_{s\rightarrow +\infty} \zeta^\star\left( k_1,\cdots,k_{r-1},k_r+1,\{1\}^s    \right),
  \]
  it follows that 
  \[
  \bigcup_{\substack{k_1\geq 2, k_2,\cdots, k_t\geq 1\\ (k_1,k_2, \cdots, k_t)\neq (2,\{1\}^{t-1})   }}  Z_{k_1,\cdots,k_t}=\left(1, \zeta^\star(2,\{1\}^{t-1})\right] ,
  \]
  \[
   \bigcup_{{k_{t+1}, k_{t+2},\cdots, k_{t+s} \geq 1  }}  Z_{k_1,\cdots,k_t, k_{t+1},\cdots, k_{t+s}}= Z_{k_1,\cdots,k_t}.  \]     
   
    \begin{lem}\label{z}
  For $k_1\geq 2, k_2,\cdots,k_r\geq 1$ and 
 $ (k_1,k_2,\cdots,k_r)\neq (2,\{1\}^{r-1})$,
  \[
  m\left(  Z_{k_1,\cdots,k_r}  \right)=\sum_{n_1\geq\cdots\geq n_{r-1}\geq n_r\geq 2} \frac{n_r-1}{n_1^{k_1}\cdots n_{r-1}^{k_{r-1}}n_r^{k_r}}.
  \]
  \end{lem}
   \noindent{\bf Proof:}
If $k_r\geq 2$, then
\[
\begin{split}
&\;\;\;\;    m\left(  Z_{k_1,\cdots,k_r}  \right)     \\
&= \zeta^\star(k_1,\cdots,k_{r-1},k_r-1)       -   \zeta^\star(k_1,\cdots,k_{r-1},k_r)        \\
& =    \sum_{n_1\geq\cdots\geq n_{r-1}\geq n_r\geq 1} \frac{n_r-1}{n_1^{k_1}\cdots n_{r-1}^{k_{r-1}}n_r^{k_r}}          \\
& =    \sum_{n_1\geq\cdots\geq n_{r-1}\geq n_r\geq 2} \frac{n_r-1}{n_1^{k_1}\cdots n_{r-1}^{k_{r-1}}n_r^{k_r}}          \\\end{split}
\] 
If $k_i\geq 2, k_{i+1}=\cdots =k_r=1$, then
\[
\begin{split}
&\;\;\;\;    m\left(  Z_{k_1,\cdots,k_r}  \right)     \\
&= \zeta^\star(k_1,\cdots,k_{i-1},k_i-1)       -   \zeta^\star(k_1,\cdots,k_{r-1},k_r)        \\
& = \sum_{n_1\geq \cdots \geq n_{i-1}\geq n_i\geq 2} \frac{n_i}{n_1^{k_1}\cdots n_{i-1}^{k_{i-1}} n_i^{k_i}}   -\sum_{n_1\geq\cdots\geq n_{r-1}\geq n_r\geq 1} \frac{1}{n_1^{k_1}\cdots n_{r-1}^{k_{r-1}}n_r^{k_r}}    .      \\
\end{split}
\]
By the argument in Section \ref{gea}, for $n_i\geq1$ one has 
\[
n_i=\sum_{n_i\geq n_{i+1}\geq \cdots \geq n_{r-1}\geq n_r\geq 1}\frac{1}{n_{i+1}\cdots n_{r-1}}.
\]
Thus for $k_i\geq 2, k_{i+1}=\cdots =k_r=1$, 
\[
\begin{split}
&\;\;\;\;    m\left(  Z_{k_1,\cdots,k_r}  \right)     \\
& = \sum_{n_1\geq \cdots \geq n_{i-1}\geq n_i\geq 2} \frac{n_i}{n_1^{k_1}\cdots n_{i-1}^{k_{i-1}} n_i^{k_i}}   -\sum_{n_1\geq\cdots\geq n_{r-1}\geq n_r\geq 1} \frac{1}{n_1^{k_1}\cdots n_{r-1}^{k_{r-1}}n_r^{k_r}}     \\
&=  \sum_{n_1\geq \cdots \geq n_{r-1}\geq n_r\geq 2} \frac{1}{n_1^{k_1}\cdots n_{i-1}^{k_{i-1}} n_i^{k_i}n_{i+1}\cdots n_{r-1}}   -\sum_{n_1\geq\cdots\geq n_{r-1}\geq n_r\geq 1} \frac{1}{n_1^{k_1}\cdots n_{r-1}^{k_{r-1}}n_r^{k_r}}       \\
& =    \sum_{n_1\geq\cdots\geq n_{r-1}\geq n_r\geq 2} \frac{n_r-1}{n_1^{k_1}\cdots n_{r-1}^{k_{r-1}}n_r^{k_r}}  .        \\\end{split}
\]

$\hfill\Box$\\

    \begin{lem}   \label{2pc}   For  $h:\mathbb{N}\rightarrow [1,+\infty)$, if the series
 \[
 \sum_{r\geq 1}\frac{1}{2^{h(r)}}
 \]
 is divergent, then  for $t\geq 2$
        \[
        \begin{split}
        &\lim_{p\rightarrow +\infty} \sum_{n_1\geq \cdots\geq n_{p-1}\geq n_p\geq  2}\frac{1}{n_1(n_1-1)(n_2-1)\cdots (n_t-1)} \cdot \frac{1}{n_t}   \\
& \cdot \left[ \frac{ 1-\left( \frac{1}{n_{t+1}}\right)^{h(t)}}{n_{t+1}-1}  \right]\cdots \left[ \frac{ 1-\left( \frac{1}{n_{p-1}}\right)^{h(p-2)}}{n_{p-1}-1}  \right] \left[ 1-\left(  \frac{1}{n_p}   \right)^{h(p-1)}    \right]  \\
&=0.
\end{split}        
        \]
        \end{lem}
       \noindent{\bf Proof: }  
        For $t+1\leq p-1$, we have 
          \[
        \begin{split}
        &\;\;\;\; \sum_{n_1\geq \cdots\geq n_{p-1}\geq n_p\geq  2}\frac{1}{n_1(n_1-1)(n_2-1)\cdots (n_t-1)} \cdot \frac{1}{n_t}   \\
&\;\;\;\; \cdot \left[ \frac{ 1-\left( \frac{1}{n_{t+1}}\right)^{h(t)}}{n_{t+1}-1}  \right]\cdots \left[ \frac{ 1-\left( \frac{1}{n_{p-1}}\right)^{h(p-2)}}{n_{p-1}-1}  \right] \left[ 1-\left(  \frac{1}{n_p}   \right)^{h(p-1)}    \right]  \\
&=\Bigg{(} \sum_{n_1= \cdots= n_{p-1}= n_p=  2} + \sum_{\substack{n_1\geq  3\\ n_2 =\cdots= n_{p-1}= n_p=  2}}+\cdots+ \sum_{\substack{n_1\geq \cdots \geq n_t\geq  3\\ n_{t+1} =\cdots= n_p=  2}} \\
&\;\;\;\;+\sum_{\substack{n_1\geq \cdots \geq n_{t+1}\geq  3\\ n_{t+2} =\cdots= n_p=  2}}  +\cdots+\sum_{n_1\geq \cdots\geq n_{p-1}\geq n_p\geq  3} \Bigg{)}   \frac{1}{n_1(n_1-1)(n_2-1)\cdots (n_t-1)} \cdot \frac{1}{n_t}     \\
&\;\;\;\; \cdot   \left[ \frac{ 1-\left( \frac{1}{n_{t+1}}\right)^{h(t)}}{n_{t+1}-1}  \right]\cdots \left[ \frac{ 1-\left( \frac{1}{n_{p-1}}\right)^{h(p-2)}}{n_{p-1}-1}  \right] \cdot \left[ 1-\left(  \frac{1}{n_p}   \right)^{h(p-1)}    \right]  \\
\end{split}
\]
\[
\begin{split}
&= \left(1-\frac{1}{2^{h(t)}} \right) \cdots \left(   1-\frac{1}{2^{h(p-2)}}     \right)  \left(   1-\frac{1}{2^{h(p-1)}}     \right) \Bigg{[}  \frac{1}{4}+ \frac{1}{2}\sum_{n_1\geq 3} \frac{1}{n_1(n_1-1)}+  \cdots\\
&\;\;\;\;+ \frac{1}{2}\sum_{n_1\geq \cdots n_{t-1}\geq 3} \frac{1}{n_1(n_1-1)(n_2-1)\cdots (n_{t-1}-1)}  \\
&\;\;\;\;+\sum_{n_1\geq \cdots\geq n_{t-1}\geq n_t\geq  3}\frac{1}{n_1(n_1-1)(n_2-1)\cdots (n_t-1)} \cdot \frac{1}{n_t}    \Bigg{]}
 \\
 &\;\;\;\; +\Bigg{(}\sum_{\substack{n_1\geq \cdots \geq n_{t+1}\geq  3\\ n_{t+2} =\cdots= n_p=  2}}  +\cdots+\sum_{n_1\geq \cdots\geq n_{p-1}\geq n_p\geq  3} \Bigg{)}   \frac{1}{n_1(n_1-1)(n_2-1)\cdots (n_t-1)} \cdot \frac{1}{n_t}     \\
&\;\;\;\; \cdot   \left[ \frac{ 1-\left( \frac{1}{n_{t+1}}\right)^{h(t)}}{n_{t+1}-1}  \right]\cdots \left[ \frac{ 1-\left( \frac{1}{n_{p-1}}\right)^{h(p-2)}}{n_{p-1}-1}  \right] \cdot \left[ 1-\left(  \frac{1}{n_p}   \right)^{h(p-1)}    \right] . \\
 \end{split}        
        \]
        As 
        \[
        \lim_{p\rightarrow+\infty}  \left(1-\frac{1}{2^{h(t)}} \right) \cdots \left(   1-\frac{1}{2^{h(p-2)}}     \right)  \left(   1-\frac{1}{2^{h(p-1)}}     \right)=0,        \]
        it suffices to show that 
        \[
        \begin{split}
        &\lim_{p\rightarrow +\infty} \Bigg{(}\sum_{\substack{n_1\geq \cdots \geq n_{t+1}\geq  3\\ n_{t+2} =\cdots= n_p=  2}}  +\cdots+\sum_{n_1\geq \cdots\geq n_{p-1}\geq n_p\geq  3} \Bigg{)}   \frac{1}{n_1(n_1-1)(n_2-1)\cdots (n_t-1)} \cdot \frac{1}{n_t}     \\
&\;\;\;\; \cdot   \left[ \frac{ 1-\left( \frac{1}{n_{t+1}}\right)^{h(t)}}{n_{t+1}-1}  \right]\cdots \left[ \frac{ 1-\left( \frac{1}{n_{p-1}}\right)^{h(p-2)}}{n_{p-1}-1}  \right] \cdot \left[ 1-\left(  \frac{1}{n_p}   \right)^{h(p-1)}    \right] =0.\\
 \end{split}        
        \]
        For $t+2\leq q \leq p-2$, we have 
           \[
        \begin{split}
        &\;\;\;\;\Bigg{(}\sum_{\substack{n_1\geq \cdots \geq n_{t+1}\geq  3\\ n_{t+2} =\cdots= n_p=  2}}  +\cdots+\sum_{n_1\geq \cdots\geq n_{p-1}\geq n_p\geq  3} \Bigg{)}   \frac{1}{n_1(n_1-1)(n_2-1)\cdots (n_t-1)} \cdot \frac{1}{n_t}     \\
&\;\;\;\; \cdot   \left[ \frac{ 1-\left( \frac{1}{n_{t+1}}\right)^{h(t)}}{n_{t+1}-1}  \right]\cdots \left[ \frac{ 1-\left( \frac{1}{n_{p-1}}\right)^{h(p-2)}}{n_{p-1}-1}  \right] \cdot \left[ 1-\left(  \frac{1}{n_p}   \right)^{h(p-1)}    \right]\\
&=  \Bigg{(}\sum_{\substack{n_1\geq \cdots \geq n_{t+1}\geq  3\\ n_{t+2} =\cdots= n_p=  2}}  +\cdots+\sum_{\substack{n_1\geq \cdots \geq n_{q}\geq  3\\ n_{q+1} =\cdots= n_p=  2}} + \sum_{\substack{n_1\geq \cdots \geq n_{q+1}\geq  3\\ n_{q+2} =\cdots= n_p=  2}} +\cdots+\sum_{n_1\geq \cdots\geq n_{p-1}\geq n_p\geq  3} \Bigg{)}   \\
&\;\;\;\;\frac{1}{n_1(n_1-1)(n_2-1)\cdots (n_t-1)} \cdot \frac{1}{n_t}     
\cdot   \left[ \frac{ 1-\left( \frac{1}{n_{t+1}}\right)^{h(t)}}{n_{t+1}-1}  \right]\cdots \left[ \frac{ 1-\left( \frac{1}{n_{p-1}}\right)^{h(p-2)}}{n_{p-1}-1}  \right] \\
&\;\;\;\;\cdot \left[ 1-\left(  \frac{1}{n_p}   \right)^{h(p-1)}    \right]\\   
&<  \Bigg{(}\sum_{\substack{n_1\geq \cdots \geq n_{t+1}\geq  3\\ n_{t+2} =\cdots= n_p=  2}}  +\cdots+\sum_{\substack{n_1\geq \cdots \geq n_{q}\geq  3\\ n_{q+1} =\cdots= n_p=  2}} \Bigg{)}   \frac{1}{n_1(n_1-1)(n_2-1)\cdots (n_t-1)} \cdot \frac{1}{n_t}     \\
&\;\;\;\;\cdot   \left[ \frac{ 1-\left( \frac{1}{n_{t+1}}\right)^{h(t)}}{n_{t+1}-1}  \right]\cdots \left[ \frac{ 1-\left( \frac{1}{n_{p-1}}\right)^{h(p-2)}}{n_{p-1}-1}  \right] \cdot \left[ 1-\left(  \frac{1}{n_p}   \right)^{h(p-1)}    \right]\\ 
&+  \Bigg{(} \sum_{\substack{n_1\geq \cdots \geq n_{q+1}\geq  3\\ n_{q+2} =\cdots= n_p=  2}} +\cdots+\sum_{n_1\geq \cdots\geq n_{p-1}\geq n_p\geq  3} \Bigg{)}   \frac{1}{n_1(n_1-1)(n_2-1)\cdots (n_t-1)} \cdot \frac{1}{n_t}     \\
&\;\;\;\;\cdot   \frac{ 1}{(n_{t+1}-1)\cdots (n_{p-1}-1)}  \\
        &\leq  \left(1-\frac{1}{2^{h(q)}}\right)\cdots \left( 1-\frac{1}{2^{h(p-1)}}   \right)   \Bigg{(}\sum_{\substack{n_1\geq \cdots \geq n_{t+1}\geq  3\\ n_{t+2} =\cdots= n_q=  2}}  +\cdots+\sum_{\substack{n_1\geq \cdots \geq n_{q}\geq  3\\ }} \Bigg{)} \\
        &\;\;\;\;  \frac{1}{n_1(n_1-1)(n_2-1)\cdots (n_t-1)} \cdot \frac{1}{n_t}     \cdot   \left[ \frac{ 1-\left( \frac{1}{n_{t+1}}\right)^{h(t)}}{n_{t+1}-1}  \right]\cdots \left[ \frac{ 1-\left( \frac{1}{n_{q}}\right)^{h(q-1)}}{n_{q}-1}  \right] \\  
        &\;\;\;\;+  \Bigg{(} \sum_{\substack{n_1\geq \cdots \geq n_{q+1}\geq  3\\ n_{q+2} =\cdots= n_p=  2}} +\cdots+\sum_{n_1\geq \cdots\geq n_{p-1}\geq n_p\geq  3} \Bigg{)}   \frac{1}{(n_1-1)^2(n_2-1)\cdots (n_{t-1}-1) (n_t-1)^2}    \\
&\;\;\;\;\cdot   \frac{ 1}{(n_{t+1}-1)\cdots (n_{p-1}-1)}  \\
  \end{split}        
        \]
        \[
        \begin{split} &\leq  \left(1-\frac{1}{2^{h(q)}}\right)\cdots \left( 1-\frac{1}{2^{h(p-1)}}   \right)   \Bigg{(}\sum_{\substack{n_1\geq \cdots \geq n_{t+1}\geq  3\\ n_{t+2} =\cdots= n_q=  2}}  +\cdots+\sum_{\substack{n_1\geq \cdots \geq n_{q}\geq  3\\ }} \Bigg{)} \\
        &\;\;\;\;  \frac{1}{n_1(n_1-1)(n_2-1)\cdots (n_t-1)} \cdot \frac{1}{n_t}     \cdot   \left[ \frac{ 1-\left( \frac{1}{n_{t+1}}\right)^{h(t)}}{n_{t+1}-1}  \right]\cdots \left[ \frac{ 1-\left( \frac{1}{n_{q}}\right)^{h(q-1)}}{n_{q}-1}  \right] \\  
        &\;\;\;\;+  \Bigg{(} \sum_{\substack{n_1\geq \cdots \geq n_{q+1}\geq  2\\ n_{q+2} =\cdots= n_p=  1}} +\cdots+\sum_{n_1\geq \cdots\geq n_{p-1}\geq n_p\geq  2} \Bigg{)}   \frac{1}{n_1^2 n_2\cdots n_{t-1} n_t^2}    \cdot   \frac{ 1}{n_{t+1}\cdots n_{p-1}}  \\
                \end{split}
        \]
        As 
        \[
        \sum_{n_1\geq \cdots n_r\geq n_{r+1}\geq 2}\frac{1}{n_1^{k_1}\cdots n_r^{k_r}n_{r+1}^{k_{r+1}}}=\zeta^\star(k_1,\cdots,k_r,k_{r+1})-\zeta^\star(k_1,\cdots,k_r)
        \]
        and 
        \[
         \sum_{n_1\geq \cdots \geq n_r\geq 2}\frac{n_r-1}{n_1^{k_1}\cdots n_r^{k_r}}=  m\left(Z_{k_1,\cdots,k_r} \right)      \]
         for any $k_1\geq 2,k_2,\cdots,k_r\geq 1$, we have 
         \[
         \begin{split}
         &    \;\;\;\; \Bigg{(} \sum_{\substack{n_1\geq \cdots \geq n_{q+1}\geq  2\\ n_{q+2} =\cdots= n_p=  1}} +\cdots+\sum_{n_1\geq \cdots\geq n_{p-1}\geq n_p\geq  2} \Bigg{)}   \frac{1}{n_1^2 n_2\cdots n_{t-1} n_t^2}    \cdot   \frac{ 1}{n_{t+1}\cdots n_{p-1}}      \\
         &=  \sum_{n_1\geq \cdots\geq n_{q+1}\geq 2}\frac{1}{n_1^2n_2\cdots n_{t-1}n_t^2n_{t+1}\cdots n_{q+1}}+\cdots +   \sum_{n_1\geq \cdots\geq n_{p-1}\geq 2}\frac{1}{n_1^2n_2\cdots n_{t-1}n_t^2n_{t+1}\cdots n_{p-1}}  \\
         &\;\;\;\;+ \sum_{n_1\geq \cdots\geq n_{p-1}\geq n_p \geq 2}\frac{1}{n_1^2n_2\cdots n_{t-1}n_t^2n_{t+1}\cdots n_{p-1}}       \\
         &=\left( \zeta^\star(2,\{1\}^{t-2},2, \{1\}^{q-t+1})  -  \zeta^\star(2,\{1\}^{t-2},2, \{1\}^{q-t})   \right)+\cdots\\
         &\;\;\;\;+ \left( \zeta^\star(2,\{1\}^{t-2},2, \{1\}^{p-t-1})  -  \zeta^\star(2,\{1\}^{t-2},2, \{1\}^{p-t-2})   \right)\\
         &\;\;\;\; +\sum_{n_1\geq \cdots\geq n_{p-1} \geq 2}\frac{n_{p-1}-1}{n_1^2n_2\cdots n_{t-1}n_t^2n_{t+1}\cdots n_{p-1}}    \\
         &=\zeta^\star(2,\{1\}^{t-2},2, \{1\}^{p-t-1})-\zeta^\star(2,\{1\}^{t-2},2, \{1\}^{q-t})  +m\left(Z_{2,\{1\}^{t-2},2, \{1\}^{p-t-1}}   \right)          \\
         &=  \zeta^\star(2,\{1\}^{t-2},2, \{1\}^{p-t-1})-\zeta^\star(2,\{1\}^{t-2},2, \{1\}^{q-t})  \\
         &\;\;\;\;+ \zeta^\star(2,\{1\}^{t-1})-    \zeta^\star(2,\{1\}^{t-2},2, \{1\}^{p-t-1})    \\
         &=\zeta^\star(2,\{1\}^{t-1})-\zeta^\star(2,\{1\}^{t-2},2, \{1\}^{q-t})    .                      \end{split}
         \]
         
         Thus   \[
        \begin{split}
        &\lim_{p\rightarrow +\infty} \Bigg{(}\sum_{\substack{n_1\geq \cdots \geq n_{t+1}\geq  3\\ n_{t+2} =\cdots= n_p=  2}}  +\cdots+\sum_{n_1\geq \cdots\geq n_{p-1}\geq n_p\geq  3} \Bigg{)}   \frac{1}{n_1(n_1-1)(n_2-1)\cdots (n_t-1)} \cdot \frac{1}{n_t}     \\
&\;\;\;\; \cdot   \left[ \frac{ 1-\left( \frac{1}{n_{t+1}}\right)^{h(t)}}{n_{t+1}-1}  \right]\cdots \left[ \frac{ 1-\left( \frac{1}{n_{p-1}}\right)^{h(p-2)}}{n_{p-1}-1}  \right] \cdot \left[ 1-\left(  \frac{1}{n_p}   \right)^{h(p-1)}    \right] \\
&\leq   \lim_{p\rightarrow +\infty} \Bigg{\{}      \left(1-\frac{1}{2^{h(q)}}\right)\cdots \left( 1-\frac{1}{2^{h(p-1)}}   \right)   \Bigg{(}\sum_{\substack{n_1\geq \cdots \geq n_{t+1}\geq  3\\ n_{t+2} =\cdots= n_q=  2}}  +\cdots+\sum_{\substack{n_1\geq \cdots \geq n_{q}\geq  3\\ }} \Bigg{)} \\
        &\;\;\;\;  \frac{1}{n_1(n_1-1)(n_2-1)\cdots (n_t-1)} \cdot \frac{1}{n_t}     \cdot   \left[ \frac{ 1-\left( \frac{1}{n_{t+1}}\right)^{h(t)}}{n_{t+1}-1}  \right]\cdots \left[ \frac{ 1-\left( \frac{1}{n_{q}}\right)^{h(q-1)}}{n_{q}-1}  \right]           \Bigg{\}}\\
          &\;\;\;\;+\zeta^\star(2,\{1\}^{t-1})-\zeta^\star(2,\{1\}^{t-2},2, \{1\}^{q-t}) \\
          &\leq \zeta^\star(2,\{1\}^{t-1})-\zeta^\star(2,\{1\}^{t-2},2, \{1\}^{q-t}) \\            \end{split}        
        \]
        for any $q\geq t+2$.
        Since 
        \[
        \lim_{q\rightarrow+\infty} \zeta^\star(2,\{1\}^{t-2},2, \{1\}^{q-t})=   \zeta^\star(2,\{1\}^{t-1})    , \]
        it follows that
            \[
        \begin{split}
        &\lim_{p\rightarrow +\infty} \sum_{n_1\geq \cdots\geq n_{p-1}\geq n_p\geq  2}\frac{1}{n_1(n_1-1)(n_2-1)\cdots (n_t-1)} \cdot \frac{1}{n_t}   \\
& \cdot \left[ \frac{ 1-\left( \frac{1}{n_{t+1}}\right)^{h(t)}}{n_{t+1}-1}  \right]\cdots \left[ \frac{ 1-\left( \frac{1}{n_{p-1}}\right)^{h(p-2)}}{n_{p-1}-1}  \right] \left[ 1-\left(  \frac{1}{n_p}   \right)^{h(p-1)}    \right]  \\
&=0
\end{split}        
        \]
for      $t\geq 2$.

          $\hfill\Box$\\

  \begin{lem}\label{div}
For $\alpha>1 $ and 
 \[
 \alpha=\eta\left( (k_1,k_2,\cdots,k_r,\cdots)    \right)=\lim_{r\rightarrow +\infty}\zeta^\star (k_1,k_2,\cdots,k_r),
 \]
 define $k_{r}(\alpha)$ as   $$k_{r}(\alpha) =k_{r}$$ for any $r\geq 1$. For $r\geq 1$ and $f:\mathbb{N}\rightarrow [1,+\infty)$, if the series
 \[
 \sum_{r\geq 1}\frac{1}{2^{f(r)}}
 \]
 is divergent, then for almost all $\alpha\in (1,+\infty)$, the inequality 
 \[
 k_{r+1}(\alpha)\geq \left\lfloor f(r)\right\rfloor+2
 \]
 holds for infinitely many $r$. Here $\left\lfloor x \right\rfloor $ means the least integer function.
   \end{lem}
   \noindent{\bf Proof:}
  It is clear that
  \[
  \lfloor f(r)\rfloor \leq f(r)\leq \lfloor f(r) \rfloor+1.
  \]
  The series
 \[
 \sum_{r\geq 1}\frac{1}{2^{f(r)}}
 \]
 is convergent is equivalent to the series
 \[
  \sum_{r\geq 1}\frac{1}{2^{\lfloor f(r)\rfloor }}
 \]
 is convergent.
 Without loss of generality, now we assume that $f(r)$ is an integer for every $r\geq 1$.
 Define 
 \[
 \mathcal{P}_r=\bigg{\{}\alpha\; \bigg{|}\; \alpha>1,k_{r+1}(\alpha)\geq f(r)+2\bigg{\}}
 \]
 and 
 \[
 \mathcal{I}(f)=\bigg{\{}\alpha\; \bigg{|}\; \alpha>1,k_{r+1}(\alpha)\geq f(r)+2\; \mathrm{for\; infinitely\; many} \;r\bigg{\}},
 \]
 then 
 \[
  \mathcal{I}(f)=\bigcap_{s\geq 1} \bigcup_{r\geq s}  \mathcal{P}_r\]
  and 
   \[
 (1,+\infty)- \mathcal{I}(f)=\bigcup_{s\geq 1} \bigcap_{r\geq s}  \left( (1,+\infty)- \mathcal{P}_r\right)\]  
 As a result, 
 \[
 m\left(  (1,+\infty)-  \mathcal{I}(f)  \right)=0 \Leftrightarrow  m\left(  \bigcap_{r\geq s}  \left( (1,+\infty)- \mathcal{P}_r\right)    \right)=0,\;\forall\, s\geq 1.
 \]
 The set $  \bigcap\limits_{r\geq s}  \left( (1,+\infty)- \mathcal{P}_r\right)  $ is an unbounded subset of $(1,+\infty)$. 
 As for $s\geq 2$, 
 \[
  \bigcap\limits_{r\geq s}  \left( (1,+\infty)- \mathcal{P}_r\right) \subseteq \bigcup_{t\geq 2}\left[\bigcap_{r\geq t}  \left( (1,+\infty)- \mathcal{P}_r\right) \bigcap  \left(1,\zeta^\star(2,\{1\}^{t-1})  \right)  \right] .  \]
 Thus 
  \[
m\left(  \bigcap_{r\geq s}  \left( (1,+\infty)- \mathcal{P}_r\right)    \right)=0,\; \forall\; s\geq 2,
\]
\[
\Uparrow
\]
 \[
  m\left(  \bigcap_{r\geq t}  \left( (1,+\infty)- \mathcal{P}_r\right)  \bigcap \,(1,\zeta^\star(2,\{1\}^{t-1})  \right)=0 ,\;\forall \,t\geq 2.
 \]
 By definition, \[\alpha=\lim_{l\rightarrow +\infty}\zeta^\star (k_1,k_2,\cdots,k_l)\in   \bigcap_{r\geq t}  \left( (1,+\infty)- \mathcal{P}_r\right)  \bigcap\; (1,\zeta^\star(2,\{1\}^{t-1}) \]
 if and only if $k_{r+1}\leq f(r)+1$ for $r\geq t$ and 
 \[
 (k_1,k_2,\cdots,k_t)\neq \left( 2,\{1\}^{t-1}  \right).
 \]

   By the order structure of multiple zeta-star values, for $\alpha=\eta\left(  (l_1,l_2,\cdots,l_r,\cdots)\right)$, one has 
  \[
  \alpha\in {Z}_{k_1,\cdots,k_r} \Leftrightarrow    (l_1,\cdots,l_r)=(k_1,\cdots,k_r)\neq (2,\{1\}^{r-1}).
  \]
  Thus \[\alpha\in   \bigcap_{r\geq t}  \left( (1,+\infty)- \mathcal{P}_r\right)  \bigcap\; (1,\zeta^\star(2,\{1\}^{t-1}) \]
   is equivalent to 
   \[
   \alpha \in \bigcup_{\substack{ k_1\geq 2,k_2,\cdots,k_t\geq 1  \\  (k_1,\cdots,k_t)\neq (2, \{1\}^{t-1}) } }\bigcup_{1\leq k_{t+i}\leq f(t+i-1)+1, \;1\leq i\leq p-t} Z_{k_1,\cdots,k_t,k_{t+1},\cdots, k_p}
   \]
   for every $p\geq t+1$.
   
   In conclusion,
    \[
  m\left(  \bigcap_{r\geq t}  \left( (1,+\infty)- \mathcal{P}_r\right)  \bigcap \,(1,\zeta^\star(2,\{1\}^{t-1})  \right)=0 ,\;\forall \,t\geq 2.
 \]
 \[
 \Updownarrow
 \]
 \[
 \lim_{p\rightarrow +\infty}m \left( \bigcup_{\substack{ k_1\geq 2,k_2,\cdots,k_t\geq 1  \\  (k_1,\cdots,k_t)\neq (2, \{1\}^{t-1}) } }\bigcup_{1\leq k_{t+i}\leq f(t+i-1)+1, \;1\leq i\leq p-t} Z_{k_1,\cdots,k_t,k_{t+1},\cdots, k_p}\right)=0.
 \]

By the order structure of multiple zeta-star values, we have 
\[
\begin{split}
&\;\;\;\;   m \left( \bigcup_{\substack{ k_1\geq 2,k_2,\cdots,k_t\geq 1  \\  (k_1,\cdots,k_t)\neq (2, \{1\}^{t-1}) } }\bigcup_{1\leq k_{t+i}\leq f(t+i-1)+1, \;1\leq i\leq p-t} Z_{k_1,\cdots,k_t,k_{t+1},\cdots, k_p}\right)          \\
&=  \sum_{\substack{ k_1\geq 2,k_2,\cdots,k_t\geq 1  \\  (k_1,\cdots,k_t)\neq (2, \{1\}^{t-1}) } }\sum_{1\leq k_{t+i}\leq f(t+i-1)+1, \;1\leq i\leq p-t}   m \left(    Z_{k_1,\cdots,k_t,k_{t+1},\cdots, k_p}\right)                    \\
\end{split}
\]

Since 
\[
m\left(Z_{k_1,\cdots,k_p}\right)=\sum_{n_1\geq \cdots\geq n_{p-1}\geq n_p\geq  2} \frac{n_p-1}{n_1^{k_1}\cdots n_{p-1}^{k_{p-1}}n_p^{k_p}},
\]
for $p\geq t+1$, we have 
\[
\begin{split}
&\;\;\;\;    m \left( \bigcup_{\substack{ k_1\geq 2,k_2,\cdots,k_t\geq 1  \\  (k_1,\cdots,k_t)\neq (2, \{1\}^{t-1}) } }\bigcup_{\substack{1\leq k_{t+i}\leq f(t+i-1)+1\\ 1\leq i\leq p-t}} Z_{k_1,\cdots,k_t,k_{t+1},\cdots, k_p}\right)        \\
&       =    \sum_{\substack{ k_1\geq 2,k_2,\cdots,k_t\geq 1  \\  (k_1,\cdots,k_t)\neq (2, \{1\}^{t-1}) } }\sum_{\substack{1\leq k_{t+i}\leq f(t+i-1)+1\\ 1\leq i\leq p-t}}   \sum_{n_1\geq \cdots\geq n_{p-1}\geq n_p\geq  2} \frac{n_p-1}{n_1^{k_1}\cdots n_{p-1}^{k_{p-1}}n_p^{k_p}}  \\
&= \sum_{\substack{k_1\geq 2,k_2,\cdots,k_t \geq 1      \\    (k_1,k_2,\cdots,k_t)\neq \left( 2,\{1\}^{t-1}  \right)          }     } \sum_{n_1\geq \cdots\geq n_{p-1}\geq n_p\geq  2}  \frac{1}{n_1^{k_1}\cdots n_t^{k_t}}\\
&\;\;\;\;\cdot \left[ \frac{ 1-\left( \frac{1}{n_{t+1}}\right)^{f(t)+1}}{ n_{t+1}-1}  \right]\cdots \left[ \frac{ 1-\left( \frac{1}{n_{p-1}}\right)^{f(p-2)+1}}{ n_{p-1}-1}  \right] \left[ 1-\left(  \frac{1}{n_p}   \right)^{f(p-1)+1}    \right]   \\
&<\left(\sum_{{k_1\geq 3,k_2,\cdots,k_t \geq 1       }     } + \sum_{{k_1\geq 2,k_2\geq 2,k_3,\cdots,k_t \geq 1       }     }  +\cdots     +  \sum_{{k_1\geq 2,k_2,\cdots,k_{t-1} \geq 1 ,k_t\geq 2      }     }    \right) \sum_{n_1\geq \cdots\geq n_{p-1}\geq n_p\geq  2}\\
& \;\;\;\; \frac{1}{n_1^{k_1}\cdots n_t^{k_t}}\cdot \left[ \frac{ 1-\left( \frac{1}{n_{t+1}}\right)^{f(t)+1}}{ n_{t+1}-1}  \right]\cdots \left[ \frac{ 1-\left( \frac{1}{n_{p-1}}\right)^{f(p-2)+1}}{n_{p-1}-1}  \right] \left[ 1-\left(  \frac{1}{n_p}   \right)^{f(p-1)+1}    \right]  \\\end{split}
\]
\[
\begin{split}
&\leq  \sum_{n_1\geq \cdots\geq n_{p-1}\geq n_p\geq  2}\frac{1}{n_1(n_1-1)(n_2-1)\cdots (n_t-1)}\left( \frac{1}{n_1}+\frac{1}{n_2}+\cdots +\frac{1}{n_t}   \right) \\
&\;\;\;\; \cdot \left[ \frac{ 1-\left( \frac{1}{n_{t+1}}\right)^{f(t)+1}}{ n_{t+1}-1}  \right]\cdots \left[ \frac{ 1-\left( \frac{1}{n_{r-1}}\right)^{f(p-2)+1}}{ n_{p-1}-1}  \right] \left[ 1-\left(  \frac{1}{n_p}   \right)^{f(p-1)+1}    \right]  \\
&\leq  \sum_{n_1\geq \cdots\geq n_{p-1}\geq n_p\geq  2}\frac{1}{n_1(n_1-1)(n_2-1)\cdots (n_t-1)} \cdot \frac{t}{n_t}   \\
&\;\;\;\; \cdot \left[ \frac{ 1-\left( \frac{1}{n_{t+1}}\right)^{f(t)+1}}{n_{t+1}-1}  \right]\cdots \left[ \frac{ 1-\left( \frac{1}{n_{p-1}}\right)^{f(p-2)+1}}{n_{p-1}-1}  \right] \left[ 1-\left(  \frac{1}{n_p}   \right)^{f(p-1)+1}    \right]  \\
\end{split}
\]
By Lemma \ref{2pc}, one has 
 \[
 \lim_{p\rightarrow +\infty}m \left( \bigcup_{\substack{ k_1\geq 2,k_2,\cdots,k_t\geq 1  \\  (k_1,\cdots,k_t)\neq (2, \{1\}^{t-1}) } }\bigcup_{1\leq k_{t+i}\leq f(t+i-1)+1, \;1\leq i\leq p-t} Z_{k_1,\cdots,k_t,k_{t+1},\cdots, k_p}\right)=0.
 \]
 $\hfill\Box$\\  
  
   \noindent{\bf  Proof of Theorem \ref{eq01} $(ii)$:}
 By the order structure of multiple zeta-star values, it is clear that 
 \[
 k_{r+1}(\alpha)\geq \left\lfloor f(r)\right\rfloor +2\]
 \[\Updownarrow \]
 \[
 \alpha\in \left( \zeta^\star\left(k_1(\alpha),k_2(\alpha),\cdots, k_r(\alpha)\right),    \zeta^\star\left(k_1(\alpha),k_2(\alpha),\cdots, k_r(\alpha),\lfloor f(r)\rfloor+1\right)   \right]
 \]
 \[
 \Downarrow
 \]
 \[
 \begin{split}
&\;\;\;\; \bigg{|} \alpha-\zeta^\star(k_1(\alpha),\cdots, k_r(\alpha))\bigg{|}\\
&\leq \sum_{n_1\geq \cdots\geq n_{r+1}\geq 2}\frac{1}{ n_1^{k_1(\alpha)}\cdots n_r^{k_r(\alpha)} n_{r+1}^{ \lfloor f(r)\rfloor+1}}\\
&<\sum_{n_1\geq  \cdots\geq n_{r+1}\geq 2}\frac{1}{ n_1^{k_1(\alpha)}\cdots n_r^{k_r(\alpha)} n_{r+1}^{  f(r)}}.
\end{split} \]
 By Lemma \ref{div}, if the series  \[
 \sum_{r\geq 1}\frac{1}{2^{f(r)}}
 \]
 is divergent, then for almost all $\alpha\in (1,+\infty)$, 
  \[
 \bigg{|} \alpha-\zeta^\star(k_1,\cdots, k_r)\bigg{|}<\sum_{n_1\geq  \cdots\geq n_{r+1}\geq 2}\frac{1}{ n_1^{k_1}\cdots n_r^{k_r} n_{r+1}^{  f(r)}}
  \]  
  for infinitely many $r,k_1,\cdots,k_r$ with $r\rightarrow +\infty$. So  \[m\left((1,+\infty)-\bigcap_{s\geq 1} \bigcup_{r\geq s} \mathcal{V}_r(f)  \right)=0.\]
   $\hfill\Box$\\  
   
   \section{Sums of  multiple series}
   
   In this section we will study the following kinds of  multiple series 
                    \[
                    \sum_{n_1\geq n_2\geq \cdots \geq n_r\geq m}\frac{1}{n_1^{k_1}n_2^{k_2}\cdots n_r^{k_r}}             \]
                    for $k_1\geq 2,k_2,\cdots,k_r\geq 1$ and $m\geq 2$.
    \begin{lem}\label{111}
    For $u>0$, one has 
    \[
    \sum_{n_1\geq n_2\geq \cdots \geq n_r\geq 1}\frac{1}{(n_1+u-1)(n_1+u)(n_2+u)\cdots (n_r+u)       }=\frac{1}{u}.
    \]
    \end{lem}
       \noindent{\bf Proof:}    For a fixed $n_2\geq 1$, one has 
       \[
       \sum_{n_1\geq n_2}\frac{1}{(n_1+u-1)(n_1+u)}=     \sum_{n_1\geq n_2}\left(\frac{1}{n_1+u-1} - \frac{1}{n_1+u}  \right)  =\frac{1}{n_2+u-1}.  \]
       By induction,
       \[
       \begin{split}
       &\;\;\;\;    \sum_{n_1\geq n_2\geq \cdots \geq n_r\geq 1}\frac{1}{(n_1+u-1)(n_1+u)(n_2+u)\cdots (n_r+u)       }    \\
       &=   \sum_{n_2\geq n_3\geq \cdots \geq n_r\geq 1}\frac{1}{(n_2+u-1)(n_2+u)(n_3+u)\cdots (n_r+u)       }         \\
       &\cdots \;\;\\\
       &=\sum_{n_r\geq 1}\frac{1}{(n_r+u-1)(n_r+u)}\\
       &=\frac{1}{u}.
       \end{split}
       \]
          $\hfill\Box$\\  
          
          Lemma \ref{111} can be viewed as a twisted version of  multiple star analogue of the Hurwitz zeta function, which are discussed in \cite{kom}.      
              
       By using the operator $\left(  \frac{d}{du} \right)^k$  on the formula in Lemma \ref{111}, we have 
       \[
       \begin{split}
       &\sum_{\substack{k_0+k_1+\cdots +k_r=k\\ k_0,k_1,\cdots,k_r\geq 0}           } \frac{k!}{k_0!k_1!\cdots k_r!} \left[ \left(  \frac{d}{du} \right)^{k_0 }      \frac{1}{n_1+u-1}\right] \left[ \left(  \frac{d}{du} \right)^{k_1 }      \frac{1}{n_1+u}\right]\cdots \left[ \left(  \frac{d}{du} \right)^{k_r }      \frac{1}{n_r+u}\right]   \\
       & = (-1)^k\frac{k!}{u^{k+1}}.   \\
       \end{split}    \]
       By simplifications here, we have 
       \[
\sum_{\substack{  k_0+k_1+\cdots+k_r=K      \\  k_0,k_1,\cdots,k_r\geq 1      }} \sum_{n_1\geq \cdots n_r\geq 1}\frac{1}{(n_1+u-1)^{k_0} (n_1+u)^{k_1}(n_2+u)^{k_2}\cdots (n_r+u)^{k_r}}=\frac{1}{u^{K-r}};
\]
As a result, Theorem \ref{m}, $(i)$ is proved. By letting $u=m-1$, Theorem \ref{m},$(ii)$ follows immediately from Theorem \ref{m}, $(i)$.

\section*{Acknowledgements}       
   The author wants to thank Sam Chow, Hidekazu Furusho,   Li Lai, Sijie Luo,  Qingchun Tian,  Ye Tian, Shengyou Wen, Jun Wu, Yufeng Wu, Bingyong Xie, Bin Zhang and Wei Zhang for helpful comments.

\end{document}